\documentclass[american,english]{article}

\usepackage[T1]{fontenc}
\usepackage[latin9]{inputenc}
\usepackage{geometry}
\usepackage{color}
\usepackage{array}
\usepackage{bm}
\usepackage{multirow}
\usepackage{amsmath}
\usepackage{amssymb}
\usepackage{graphicx}
\usepackage{epstopdf}
\usepackage{setspace}
\usepackage[numbers]{natbib}
\PassOptionsToPackage{normalem}{ulem}
\usepackage{ulem}
\usepackage{xcolor}
\usepackage{relsize}
\usepackage[english]{babel}

\usepackage{subfig}
		
\providecommand{\tabularnewline}{\\}

\usepackage{indentfirst}
\usepackage{authblk}

\author[1]{Shanglong Zhang}
\author[2]{Arun L. Gain}
\author[1]{Juli{\'a}n A. Norato\thanks{Corresponding author: \tt julian.norato@uconn.edu}}

\affil[1]{Department of Mechanical Engineering, University of Connecticut}
\affil[2]{VPDS Inc.}

\title{Adaptive Mesh Refinement for Topology Optimization with Discrete Geometric
	Components}

\date{}

\begin{document}

\maketitle 

\begin{abstract}
This work introduces an Adaptive Mesh Refinement (AMR) strategy for
the topology optimization of structures made of discrete geometric
components using the geometry projection method. Practical structures
made of geometric shapes such as bars and plates typically exhibit
low volume fractions with respect to the volume of the design region
they occupy. To maintain an accurate analysis and to ensure well-defined
sensitivities in the geometry projection, it is required that the
element size is smaller than the smallest dimension of each component.
For low-volume-fraction structures, this leads to finite element meshes
with very large numbers of elements. To improve the efficiency of
the analysis and optimization, we propose a strategy to adaptively
refine the mesh and reduce the number of elements by having a finer
mesh on the geometric components, and a coarser mesh away from them.
The refinement indicator stems very naturally from the geometry projection
and is thus straightforward to implement. We demonstrate the effectiveness
of the proposed AMR method by performing topology optimization for
the design of minimum-compliance and stress-constrained structures
made of bars and plates.
\end{abstract}


\section{Introduction\label{sec:Introduction}}

Topology optimization is a powerful tool to explore structural designs
given performance and resource requirements. By relating the design
parameterization with a fixed discretization for analysis, efficient
optimization algorithms can be used to determine the optimal material
distribution. The prevalent classes of methods for topology optimization
are the density-based method and the level-set method (cf., \cite{sigmund2013topology}).
In the former, the design region is discretized into a voxel-like
grid using the finite element mesh, with a continuous pseudodensity
variable assigned to each element that indicates the presence or absence
of material at that element. In level-set methods, the boundary of
the design is represented by the zero level set of a function, and
the optimal design is obtained by evolving the zero level set. Both
classes of techniques render organic designs, which then need to be
translated into a Computer-Aided Design (CAD) model adding considerations
for manufacturability. Design modifications are inevitably introduced
in this translation, hence the final design departs from its optimum.
This departure is more pronounced when the structure is fabricated
using stock material by joining components of fixed shape but variable
dimensions. 

To address this difficulty, several works have recently been proposed
in topology optimization to produce designs made exclusively of discrete
geometric components. Since these geometric components have high-level
geometry representations, these techniques ease the transition to
a CAD model that more closely resembles the optimal topology. One
category of these techniques is the geometry projection method \cite{bell2012geometry,norato2015geometry,zhang2016geometry,deng2016design,zhang2017optimal,zhang2017stress,zhanggeometry,kazemi2018topology,norato2018topology}.
This technique smoothly maps the high-level description of the geometric
components onto a fixed density grid at each optimization iteration
to circumvent re-meshing upon design changes. This method has been
successfully applied to the design of structures made of straight
bars \cite{norato2015geometry}, flat plates \cite{zhang2016geometry,zhang2017optimal}
curved plates \cite{zhanggeometry} and geometric primitives represented
by the super-formula \cite{norato2018topology}. It has been applied
both to compliance minimization problems as well as problems with
stress considerations \cite{zhang2017stress}. To further improve
the manufacturability of the optimal design, some geometric constraints
have been introduced, such as bounds on the dimensions of the components
\cite{norato2015geometry,zhang2016geometry}, ensuring a minimum
separation between geometric components to allow for tool access \cite{zhang2017optimal},
and ensuring geometric components are fully contained withing irregular,
non-convex design regions to avoid impractical cuts \cite{zhanggeometry}.
The geometry projection method has also been used for the topology
optimization of multimaterial lattices \cite{watts2017geometric}
and structures \cite{kazemi2018topology}. 

The other major category of topology optimization techniques with
geometric components is the Moving Morphable Components (MMC) method
\cite{guo2014doing,deng2016design,guo2016explicit,zhang2016minimum,guo2017self,zhang2017structural,hoang2017topology,zhang2018topology,zhang2018moving},
along with its related technique of Moving Morphable Voids (MMV) \cite{zhang2017structural,zhang2017explicit,zhang2018amoving}.
In these techniques, the geometric solid components or voids are represented
via a topological description function, which corresponds to a level
set description of the boundaries. This function is then mapped onto
the analysis domain via a smooth Heaviside function. 

Topology optimization techniques typically employ a fixed finite element
mesh with uniform (or close-to-uniform) element size. The element
size should not be too large, otherwise the design and the structural
response cannot be accurately captured. In the case of thin structural
components, it is desirable to have multiple elements across the component
thickness. In the context of geometry projection methods, a further
requirement is imposed on the element size that requires the maximum
element size must be less than half of the thickness of the component
to ensure well-defined, continuous sensitivities \cite{norato2015geometry,zhang2016geometry}.
These requirements dictate the maximum element size that can be used
in geometry projection methods to obtain good results. Structures
made of stock material often occupy a relatively small fraction of
the volume of the design region; and the thickness of the components
is very small in relation to the dimensions of the design region.
If a uniform element size mesh is used, the foregoing mesh size requirements
together with the volume fractions and dimension ratios encountered
in practice for structures made of stock material lead to very large
meshes with tens of millions of elements. Therefore, a natural idea
is to employ a coarse mesh for regions where there is pure solid or
void, and a fine mesh near the boundaries of the geometric components.
In this work, we propose an Adaptive Mesh Refinement (AMR) strategy
for the geometry projection method to reduce the total number of elements
in the mesh while keeping the same accuracy for the geometry projection
and the analysis.

Several AMR approaches have been proposed in the literature to improve
the computational efficiency and solution accuracy in density-based
topology optimization. Most methods use some function of the density
field and/or an analysis error estimator as a refinement indicator.
The earliest of these methods is presented in \cite{Maute1995},
which presents a strategy commonly used in subsequent methods. In
this strategy, a complete optimization is performed on a coarse, uniform
size mesh, whose resulting density field is used to refine the mesh
and perform a subsequent optimization. Some works (e.g., \cite{ArantesCosta2003})
also use the mesh quality and an analysis error estimator computed
on the previous mesh as additional criteria for refinement. Other
density-based AMR techniques that follow the optimize $\to$ refine
strategy based on one or more of the aforementioned refinement criteria
include \cite{Stainko2006,Bruggi2011,Wallin2012,Wang2014,Nguyen-Xuan2017}.
All of these methods aim to have a coarser mesh in the void region;
some use a finer mesh for regions that are solid or near the boundary,
while others only refine a band of elements around the boundary, which
leads to additional computational savings.

A different AMR strategy consists of adaptively refining the mesh
upon each design change in the optimization. Such a technique is proposed
in \cite{wang2010dynamic}, where it is shown that this strategy
leads to better designs, since in the foregoing strategy the early
optimization runs on coarse meshes can lock into poor local minima.
Another contribution of this work is that an additional layer of elements
near the material boundary is marked for refinement to improve the
optimality of the final design. This approach allows de-refinement
to coarsen the fine mesh in void regions as the material boundary
propagates in the optimization. Other techniques that employ the design
step $\to$ refine strategy include \cite{Nana2016,Jensen2016,jensen2016anisotropic,Panesar2017,Troya2018}.
It is worth noting that the works in \cite{Jensen2016,Troya2018}
consider stress-based topology optimization and thus have an additional
motivation for doing AMR, which is to improve the accuracy of the
stresses.

All of the aforementioned AMR schemes are for density-based topology
optimization methods. In \cite{zhang2018amoving}, AMR is performed
at the boundary of holes described via B-splines in the MMV method
to improve the accuracy of stress calculation. Elements that are cut
by the material boundary, which is given by the boundary of the B-spline,
are marked for further refinement. The analysis is solved using the
extended finite element method (XFEM).

In this work, we propose an AMR strategy for topology optimization
with discrete geometric components using the geometry projection method.
To this end, we combine and adapt several of the aforementioned AMR
strategies employed in density-based techniques. In our method, the
mesh is refined and coarsened upon each design change in the optimization
based solely on the geometric description of the components. The geometry
projection method provides a natural way to define the refinement
indicator. We only refine elements on a band around the material boundary,
leading to a finer mesh near the material boundary and a coarser mesh
at the pure solid or void regions. We note that we do not employ any
error estimator from the analysis as refinement indicator, since our
primary goal is to reduce the total number of elements for large problems.
The remainder of the paper is organized as follows. Section 2 briefly
describes the geometry projection method. Section 3 presents the AMR
strategy. The optimization and computer implementation are detailed
in Section 4. We provide numerical examples in Section 5 to demonstrate
the effectiveness of the proposed method. Finally, we draw conclusions
in Section 6. 

\section{Geometry projection}

For the topology optimization with geometric components, we map the
geometric components onto a fixed analysis mesh to circumvent re-meshing
upon design changes. To this end, we employ the geometry projection
method \cite{norato2004geometry,norato2015geometry} to establish
a differentiable map between the parametric description of the geometry
and a density field defined over a fixed finite element grid. With
the projected density, we can perform the finite element primal and
sensitivity analyses in a manner similar to the density-based method.
Since the geometry projection is smooth, a continuous design sensitivity
of the geometric components can be readily computed via the chain
rule. 

In what follows, we briefly describe the geometry projection technique,
starting with the geometry projection of a single geometric component.
We first consider a ball-shaped sampling window $B_{\mathbf{x}}^{R}$
of radius $R$ centered at $\mathbf{x}$. The projected density at
$\mathbf{x}$ is computed as the area fraction of the circular segment
(or volume fraction of spherical cap in 3-d) of height $R-d$ (see
Fig. \ref{fig:geometry_projection}):
\begin{equation}
\rho_{i}(\mathbf{x}):=\frac{\mid B_{\mathbf{x}}^{R}\cap\omega_{i}\mid}{\mid B_{\mathbf{x}}^{R}\mid}\label{eq:geometry_projection}
\end{equation}
where $\omega_{i}$ denotes the region of space occupied by component
$i$, and $\mid\cdot\mid$ is the area (or volume). Assuming that
the portion of $\partial\omega_{i}$ that intersects the ball, $\partial\omega_{i}\cap B_{\mathbf{x}}^{R}$
, can be approximated as a straight line in 2-d, the projected density
of Eq.$\ $\ref{eq:geometry_projection} can be computed as a function
of the signed distance $d(\mathbf{x})$ from $\mathbf{x}$ to the
boundary of the component $\partial\omega_{i}$ as:
\[
\rho_{i}(d(\mathbf{x}),R):=\left\{ \begin{array}{ll}
0 & \mbox{if }d>R\\
\frac{1}{\pi R^{2}}\Big[R^{2}\arccos(\frac{d}{R})-d\sqrt{R^{2}-d^{2}}\,\Big] & \mbox{if }-R\leq d\leq R\\
1 & \mbox{if }d<-R
\end{array}\right.
\]

Conversely, if for 3-dimensional problems we assume that $\partial\omega_{i}\cap B_{\mathbf{x}}^{R}$
can be approximated as a flat surface, the projected density can be
computed as:
\[
\rho_{i}(d(\mathbf{x}),R):=\left\{ \begin{array}{ll}
0 & \mbox{if }d>R\\
\frac{1}{2}+\frac{d^{3}}{4R^{3}}-\frac{3d}{4R} & \mbox{if }-R\leq d\leq R\\
1 & \mbox{if }d<-R
\end{array}\right.
\]

The argument $\mathbf{x}$ is removed from the signed distance on
the right-hand side of the above expressions for conciseness. The
signed distance is a function the vector of design variables $\mathbf{z}_{i}$
that parameterize component $i$, i.e., \textbf{$d(\mathbf{x},\mathbf{z}_{i})$}.
Different expressions of $d$ can be derived for different component
shapes. For instance, the signed distance function of bars, flat plates,
and curved plates modeled with offset surfaces are given in \cite{norato2015geometry},
\cite{zhang2016geometry}, and \cite{zhanggeometry}; for the sake
of brevity, and since we use several of these component geometries
in this work, we refer the reader to the afore cited works for these
expressions. The projected density is computed every time the design
changes.

\begin{figure}[h]
\begin{centering}
\includegraphics[width=0.65\textwidth]{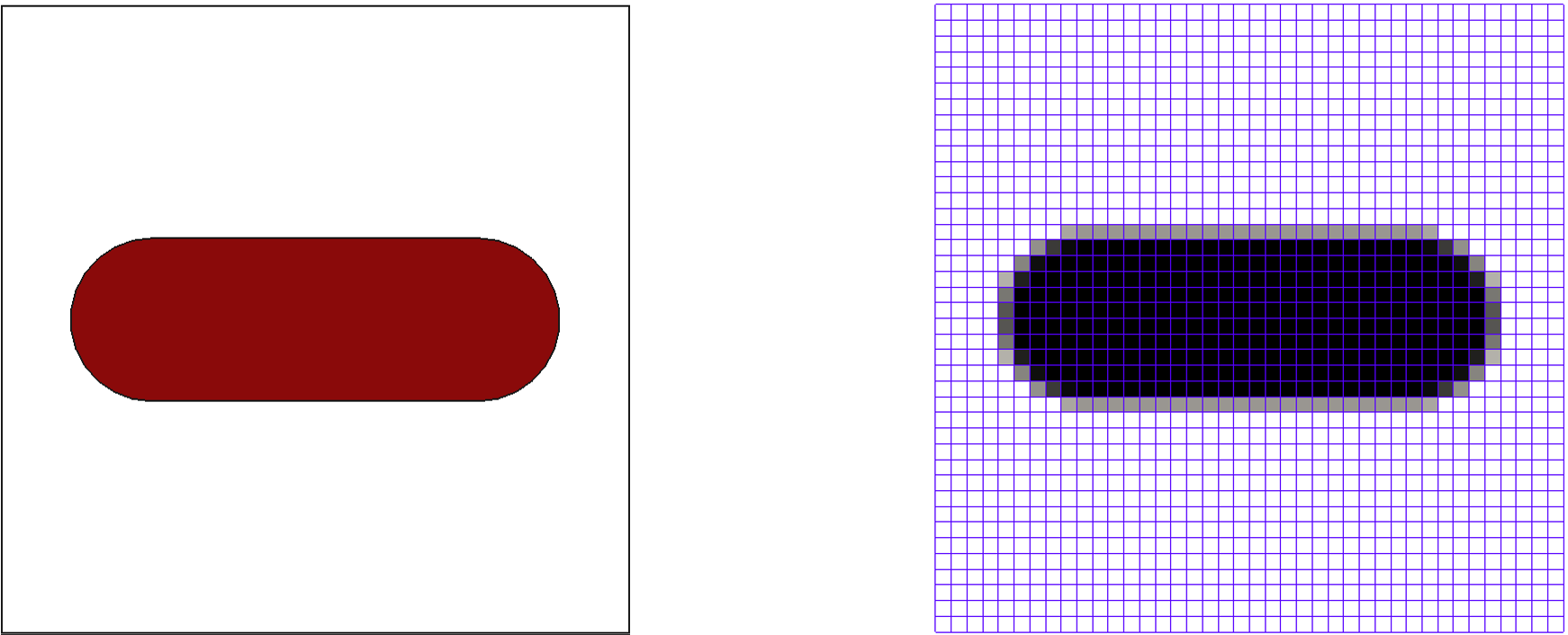}
\par\end{centering}
\caption{Geometry projection of a solid bar onto a uniform finite element grid\label{fig:bar_projection}}

\end{figure}

Fig.$\ $\ref{fig:bar_projection} shows the projected density of
a bar (modeled as an offset surface of a straight line segment, cf.,$\ $\cite{norato2015geometry}---i.e.,
the line segment is the medial axis of the bar). The idea behind the
AMR scheme we propose in this work is simple: the projected density
can be used as a refinement indicator. For instance, elements that
are fully solid or fully void ($\rho=1$ or $\rho=0$, respectively)
can be meshed with a coarse mesh, while elements with an intermediate
density ($0<\rho<1$) can be meshed with a fine mesh. We do not directly
use the projected density as a refinement indicator for reasons that
we elaborate on later---however, this is the essential idea.

%
%
%
%
%

\begin{figure}[h]
	\includegraphics[clip, trim=2cm 20cm 2cm 2cm, width=1.00\textwidth]{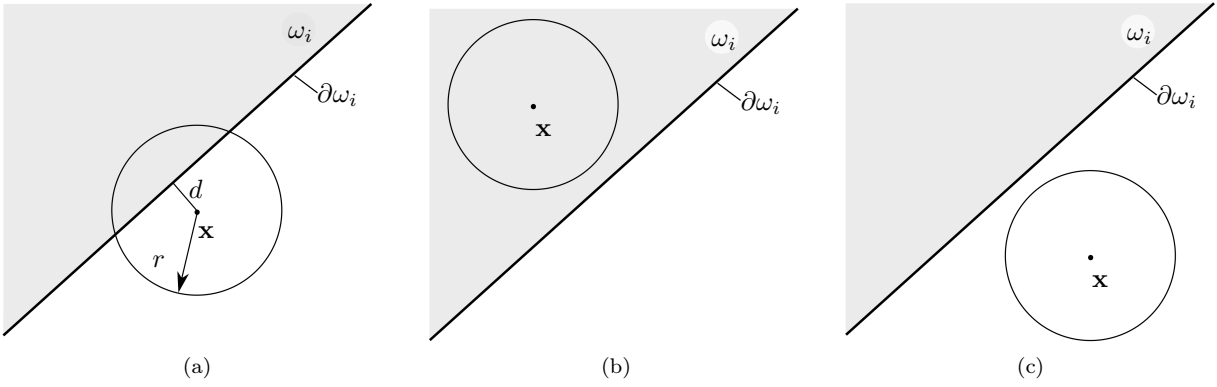}
\caption{Geometry projection for various points $\mathbf{x}$ with respect
to a geometric component $i$: (a) partially intersecting $\omega_{i}$,
with $0<\rho_{i}<1$; (b) inside of $\omega_{i}$, with $\rho_{i}=1$;
and (c) outside of $\omega_{i}$, with $\rho_{i}=0$.}
\end{figure}
As in ersatz material methods, the density is used to modified the
material properties to reflect a fully or partially solid material,
or void. In addition to the projected density of Eq.$\ $\ref{eq:geometry_projection},
we ascribe a size variable $\alpha_{i}\in[0,1]$ to each geometric
component, which is penalized in the spirit of Solid Isotropic Material
Penalization (SIMP). This is a unique and key feature of the geometry
projection method that allows the optimizer to entirely remove a geometric
component $i$ from the design. Using the projected density $\rho_{i}$
and the size variable $\alpha_{i}$, we now define an \emph{effective
density} as:
\begin{equation}
\hat{\rho}_{i}(\mathbf{x},\mathbf{z}_{i},s,R):=\alpha_{i}^{s}\rho_{i}(\mathbf{x},\mathbf{z}_{i},R)\label{eq:effective-density}
\end{equation}
where $s$ is the SIMP penalization power. Thanks to the size variable,
a component $i$ may have a projected density of unity at a point
$\mathbf{x}$, but if its size variable $\alpha_{i}$ equals zero,
the effective density of that component at $\mathbf{x}$ is also zero. 

We note that the effective density constitutes an implicit geometric
representation of the component. Therefore, when multiple components
overlap, we perform their Boolean union, which for implicit representations
corresponds to the maximum function, i.e.: 
\begin{equation}
\tilde{\rho}(\mathbf{x},\mathbf{z},s,R):=\max_{i}\hat{\rho}_{i}(\mathbf{x},\mathbf{z}_{i},s,R),\quad i=1,\cdots,n\label{eq:composite_density_max}
\end{equation}
In this expression, $\tilde{\rho}(\mathbf{x},\mathbf{z},s,R)$ is
a \emph{composite density} from all $n$ components, and $\mathbf{z}:=[\mathbf{z}_{1}^{T},\cdots,\mathbf{z}_{n}^{T}]^{T}$
denotes the vector of design variables, noting that $\mathbf{z}_{i}$
now includes $\alpha_{i}$. Since the maximum function is not differentiable,
we replace it with a smooth approximation so that we can use efficient
gradient-based optimizers. Here, we use a lower-bound Kreisselmeier--Steinhauser
(KS) function, which has the form

\begin{equation}
\max_{i}x_{i}\approx KS_{\max}(\mathbf{x}):=\frac{1}{k}\ln\left(\frac{1}{n}\sum_{i=1}^{n}e^{kx_{i}}\right)\label{eq:LKS}
\end{equation}
with aggregation coefficient $k$. We apply this approximation to
the composite density $\tilde{\rho}$ in Eq.~\ref{eq:composite_density_max}
to obtain

\begin{equation}
\tilde{\rho}_{ks}(\mathbf{x},\mathbf{z},s,R)\approx\tilde{\rho}_{KS}(\mathbf{x},\mathbf{z},s,R,k):=KS_{\max}(\hat{\bm{\rho}}(\mathbf{x},\mathbf{z},s,R))\label{eq:composite-density-ks}
\end{equation}
where $\hat{\mathbf{\bm{\rho}}}$ is the vector of effective densities.
With this smooth composite density, we now define the elasticity tensor
of the ersatz material as:

\begin{equation}
\mathbb{C}(\mathbf{x},\mathbf{z},s,R):=\left[\rho_{\min}+\tilde{\rho}_{ks}(\mathbf{x},\mathbf{z},s,R)(1-\rho_{\min})\right]\mathbb{C}_{o}\label{eq:ersatz_material}
\end{equation}
The small lower bound $0<\rho_{\min}\ll1$ prevents an ill-posed analysis.
In the finite element analysis, we use a uniform composite density
for each element, with the projected density computed at its centroid.

\section{Adaptive mesh refinement\label{sec:Adaptive-mesh-refinement}}

\begin{figure}[h]
\subfloat[\label{fig:quadtree_demo}]{\begin{centering}
\includegraphics[width=0.21\textwidth]{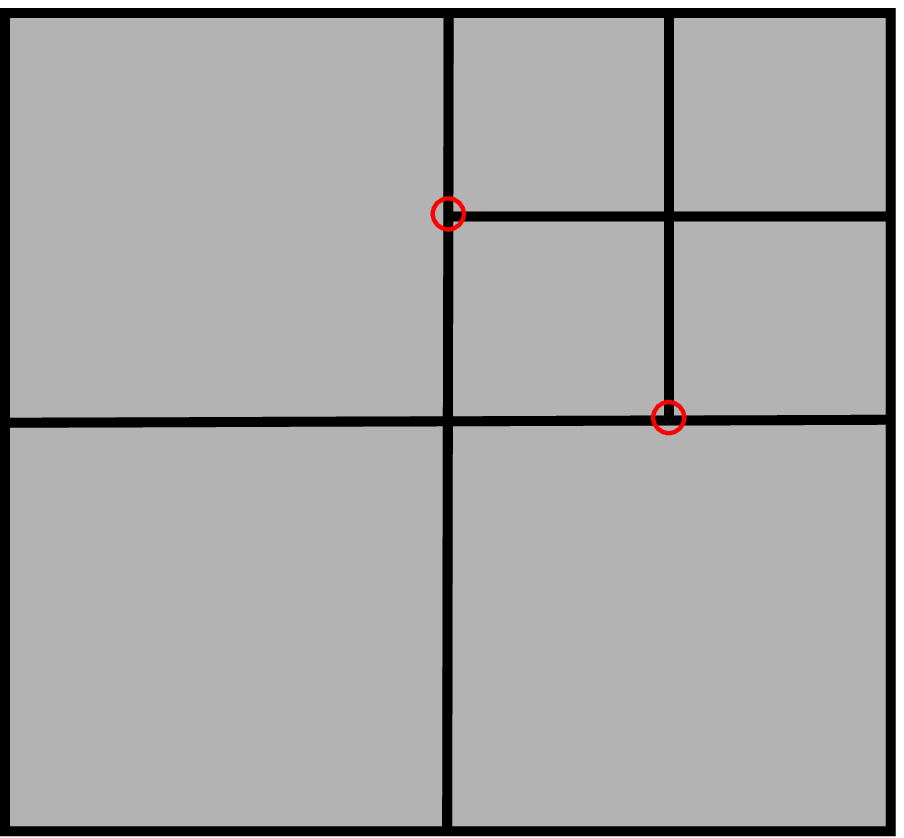}
\par\end{centering}
}\subfloat[\label{fig:octree_demo}]{\begin{centering}
\includegraphics[width=0.25\textwidth]{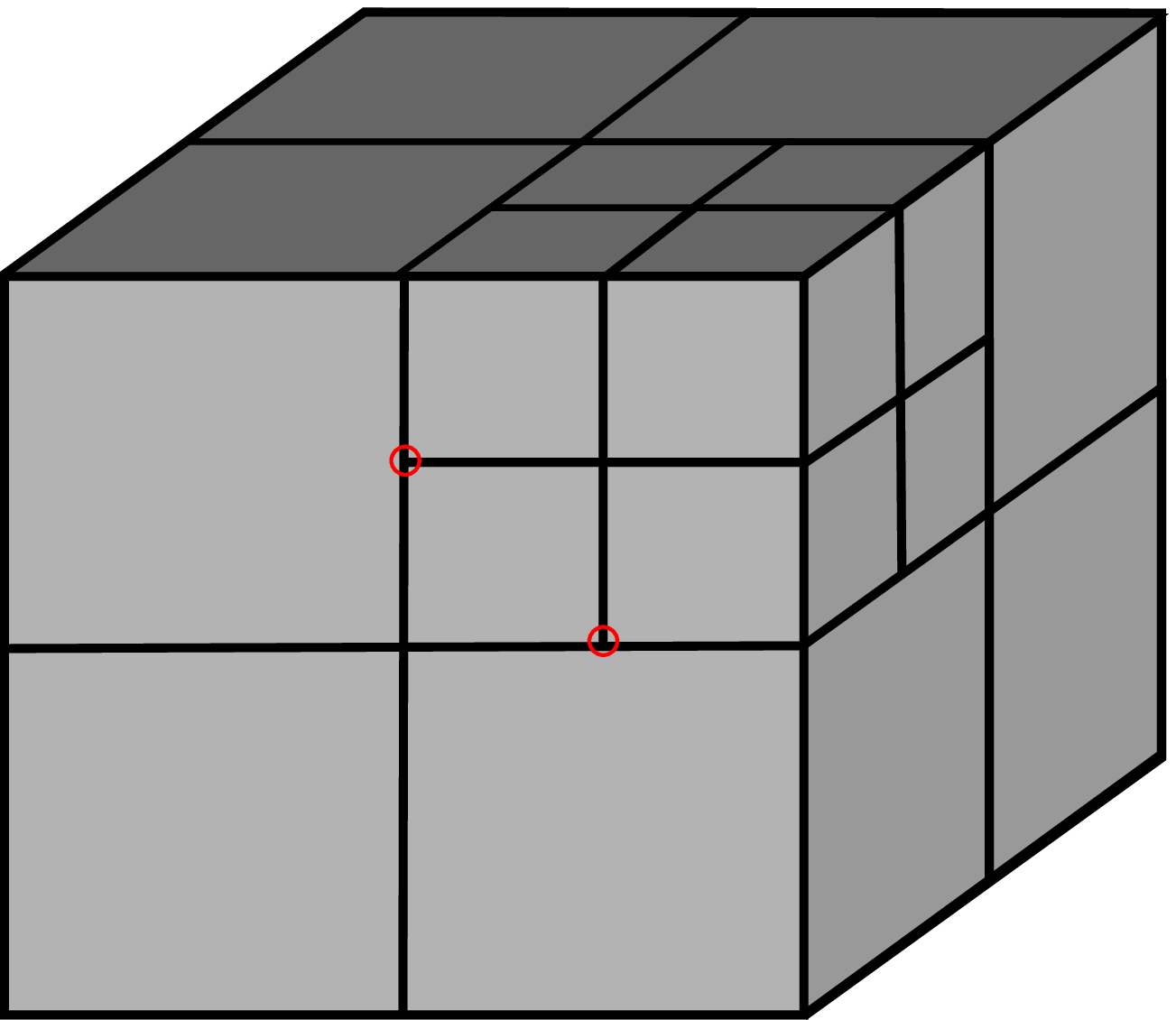}
\par\end{centering}
}\subfloat[\label{fig:quadtree_vio}]{\begin{centering}
\includegraphics[width=0.21\textwidth]{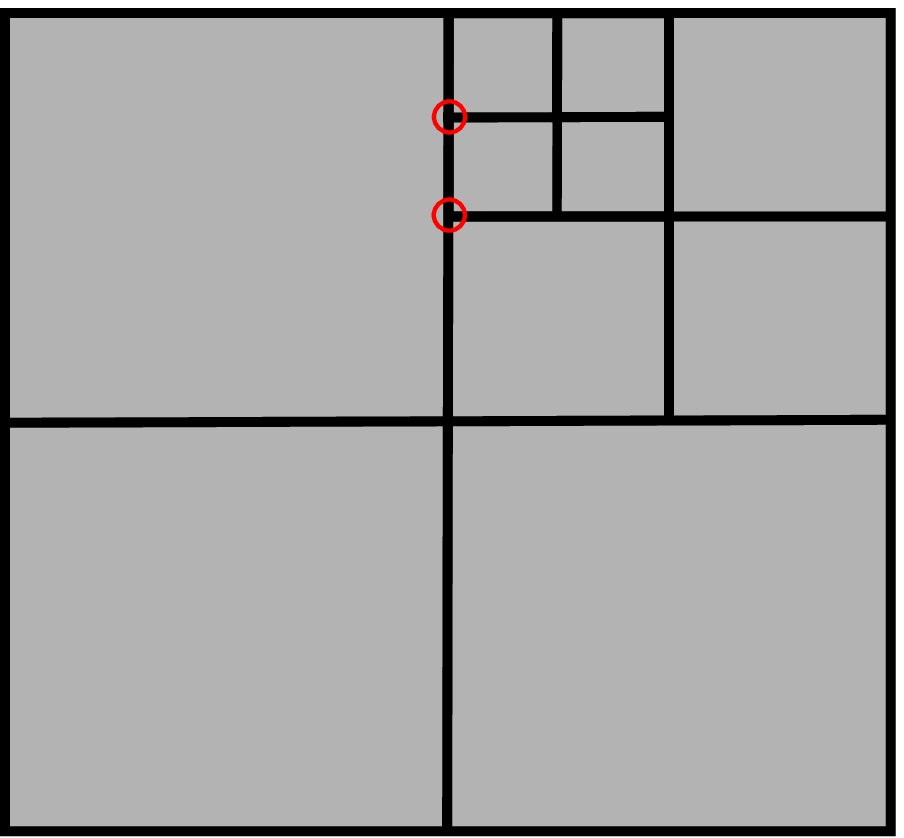}
\par\end{centering}
}\subfloat[\label{fig:octree_vio}]{\begin{centering}
\includegraphics[width=0.25\textwidth]{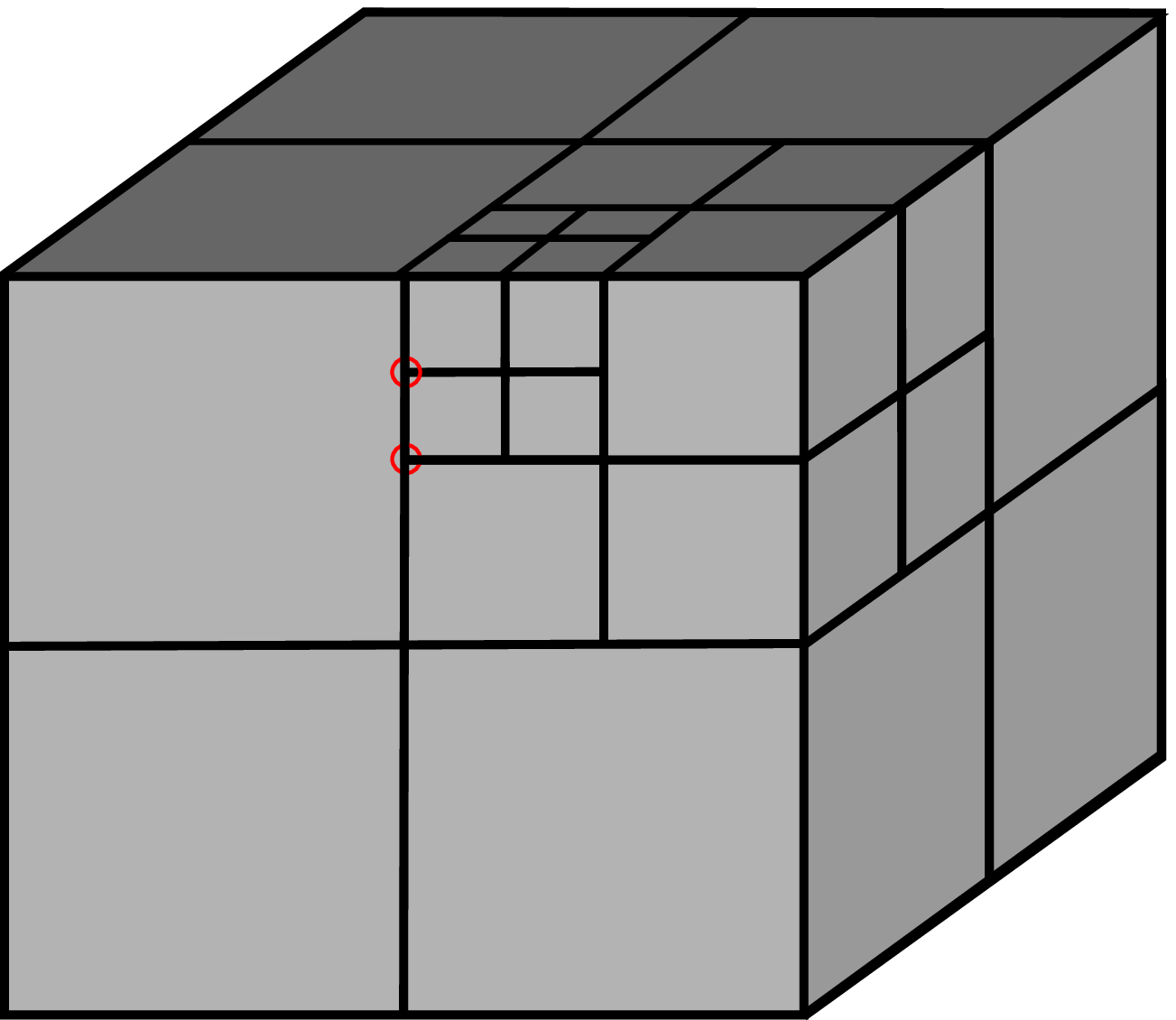}
\par\end{centering}
}

\caption{Mesh refinement by subdivision of (a and c) quadrilateral and (b and
d) hexahedral elements. The refinements in (a) and (b) satisfy the
single-level mesh-incompatibility requirement, while the refinements
in (c) and (d) violate it.\label{fig:mesh-refinement}}
\end{figure}

For a convenient demonstration of the mesh refinement strategy we
describe in the sequel, we discretize the design region using bilinear
quadrilateral elements for 2-dimensional problems and trilinear hexahedral
elements for 3-dimensional problems. The mesh refinement is performed
by successively subdividing quadrilaterals in four using a quadtree
strategy, and octahedrals in eight using an octree strategy. In Figs.\ \ref{fig:quadtree_demo}
and \ref{fig:octree_demo}, the element on the top right corner is
refined by subdividing the parent element into smaller children elements.
When neighboring elements have different levels of refinement, there
will be `hanging nodes' along element edges (indicated with circles
in Fig.\ \ref{fig:mesh-refinement}). These hanging nodes introduce
incompatibility to the finite element solution across element edges,
therefore constraint equations on the degrees of freedom must be imposed
to ensure continuity. As usual in these refinement schemes, neighboring
elements are required to differ by no more than one level of refinement
for ease of implementation. Figs\ \ref{fig:quadtree_vio} and \ref{fig:octree_vio}
show examples of mesh refinements that violate this single-level mesh-incompatibility
requirement. To satisfy this requirement, if two neighboring elements
already differ by a level of refinement and the element with the finer
mesh is marked for refinement, then the other element will be marked
for refinement too. In the case of Figs\ \ref{fig:quadtree_vio}
and \ref{fig:octree_vio}, refinement is required for the left-top
and left-top-front elements respectively.

To mark an element for refinement we require a refinement criterion.
In the case of topology optimization with discrete geometric components,
we seek a refinement indicator that accounts not only for the geometric
parameters of the component, but also for its size variable. We first
discuss for simplicity the case of a single component. In this case,
if the component has a size variable of unity, we would like the underlying
mesh to be coarse both inside and outside of the component and fine
near its boundaries; and if the size variable is zero, we would like
the underlying mesh to be coarse in the region occupied by the component.
Therefore, a suitable refinement indicator for a single component
is the effective density of Eq.$\ $\ref{eq:effective-density}. That
is, an element $e$ is marked for refinement if its effective density
satisfies 
\begin{equation}
0<\hat{\rho}^{e}\leq\rho_{th}\lesssim1\label{eq:eff-density-indicator}
\end{equation}
where $\rho_{th}$ is a threshold effective density above which the
element is considered to be close to solid. 

Based on this refinement indicator, the proposed AMR process for a
single component is as follows. 1) We start with a mesh whose elements
all have a size corresponding to the coarsest level, and we compute
the effective density for all the elements; 2) elements are marked
for refinement if they satisfy Eq.$\ $\ref{eq:eff-density-indicator};
and 3) marked elements are refined once by subdivision. If additional
levels of refinement are desired, then the process is repeated, i.e.,
the effective density is computed for all elements in the refined
mesh of step 2, and elements are subsequently marked for refinement
based on the aforementioned criterion, but also if they violate the
single-level mesh-incompatibility requirement. 

\begin{figure}[h]
\subfloat[\label{fig:coarse_grid_projection}]{\includegraphics[width=0.25\textwidth]{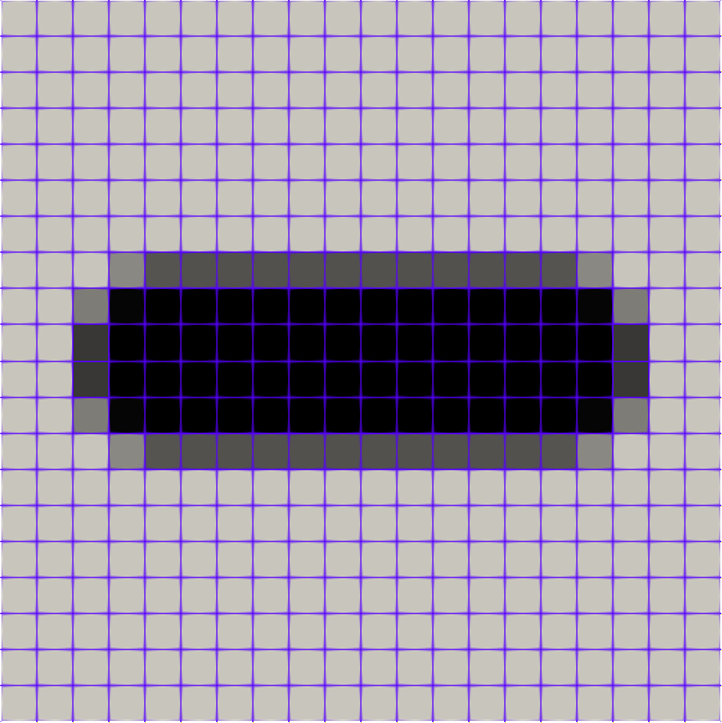}

}\subfloat[\label{fig:coarse_grid_indicator}]{\includegraphics[width=0.25\textwidth]{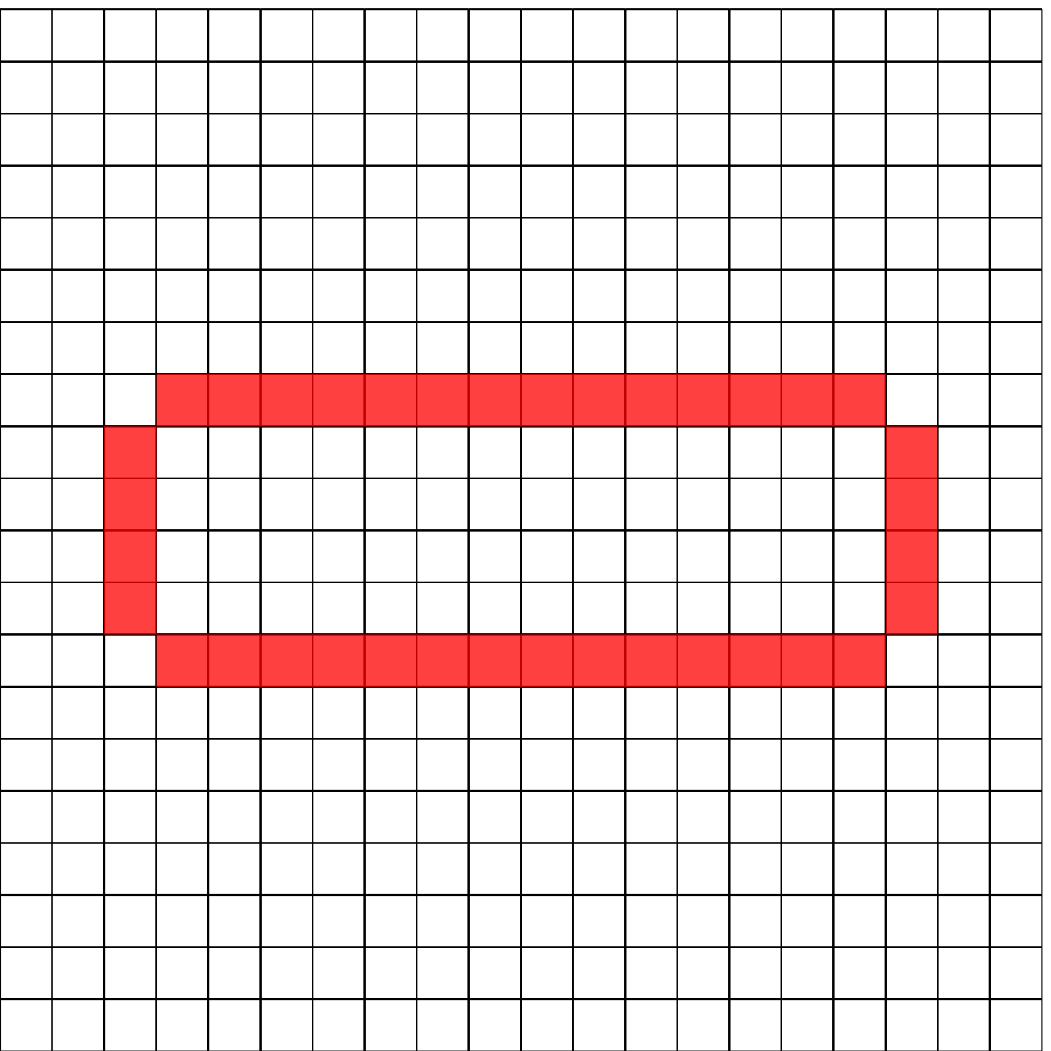}

}\subfloat[\label{fig:refined_mesh_l1}]{\includegraphics[width=0.25\textwidth]{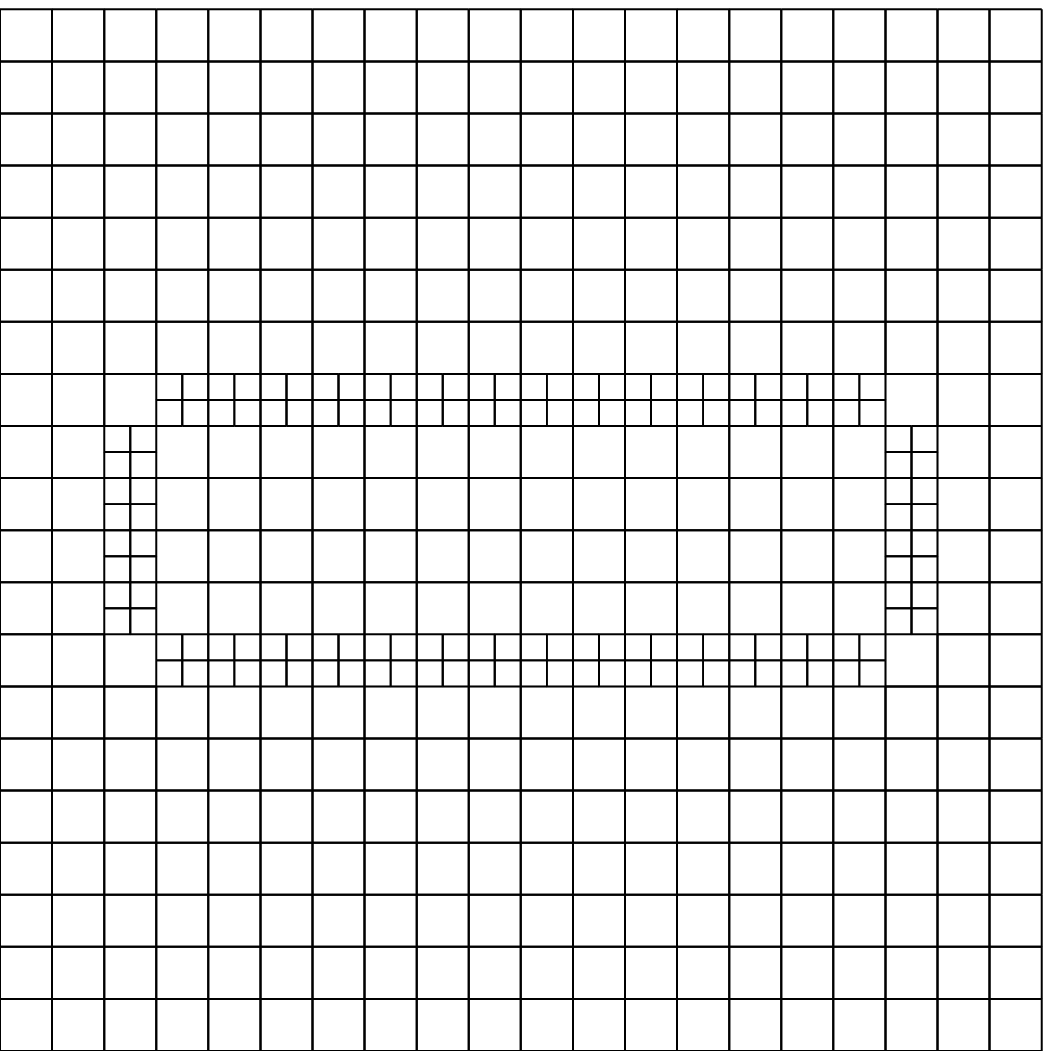}

}\subfloat[\label{fig:l1_grid_projection}]{\includegraphics[width=0.25\textwidth]{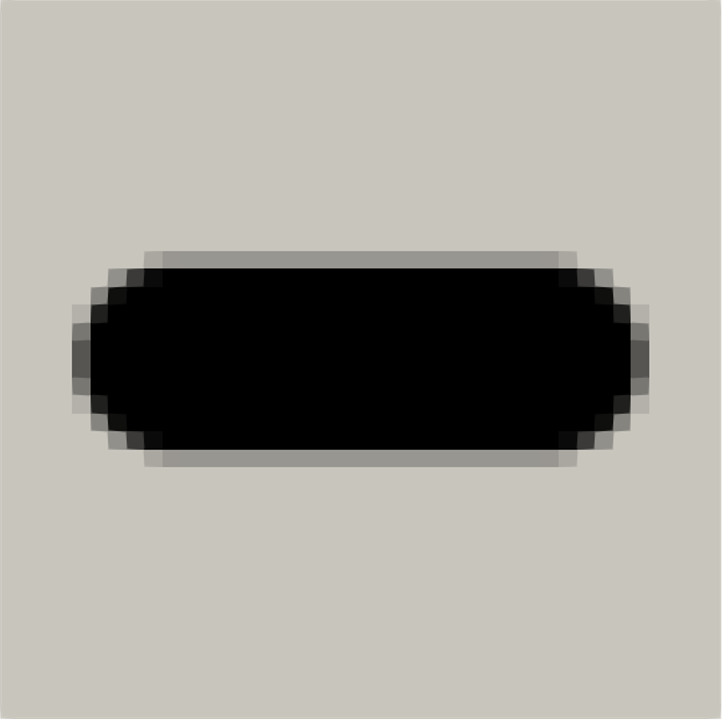}

}

\subfloat[]{\includegraphics[width=0.25\textwidth]{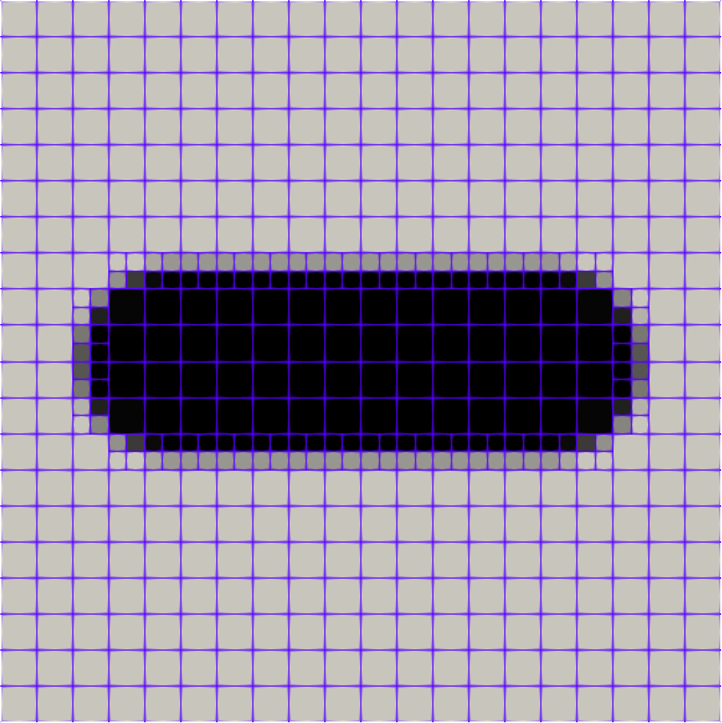}

}\subfloat[]{\includegraphics[width=0.25\textwidth]{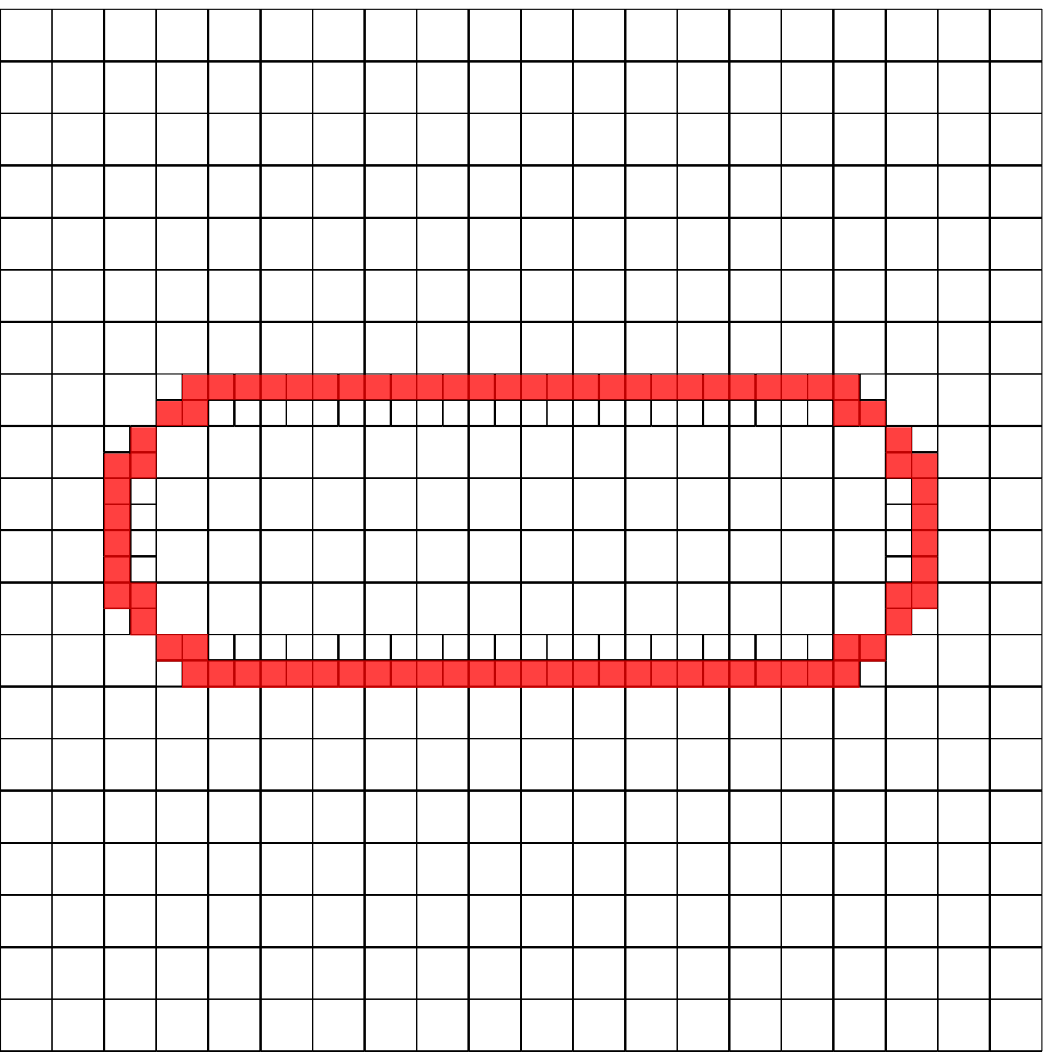}

}\subfloat[\label{fig:refine_mesh_lv2}]{\includegraphics[width=0.25\textwidth]{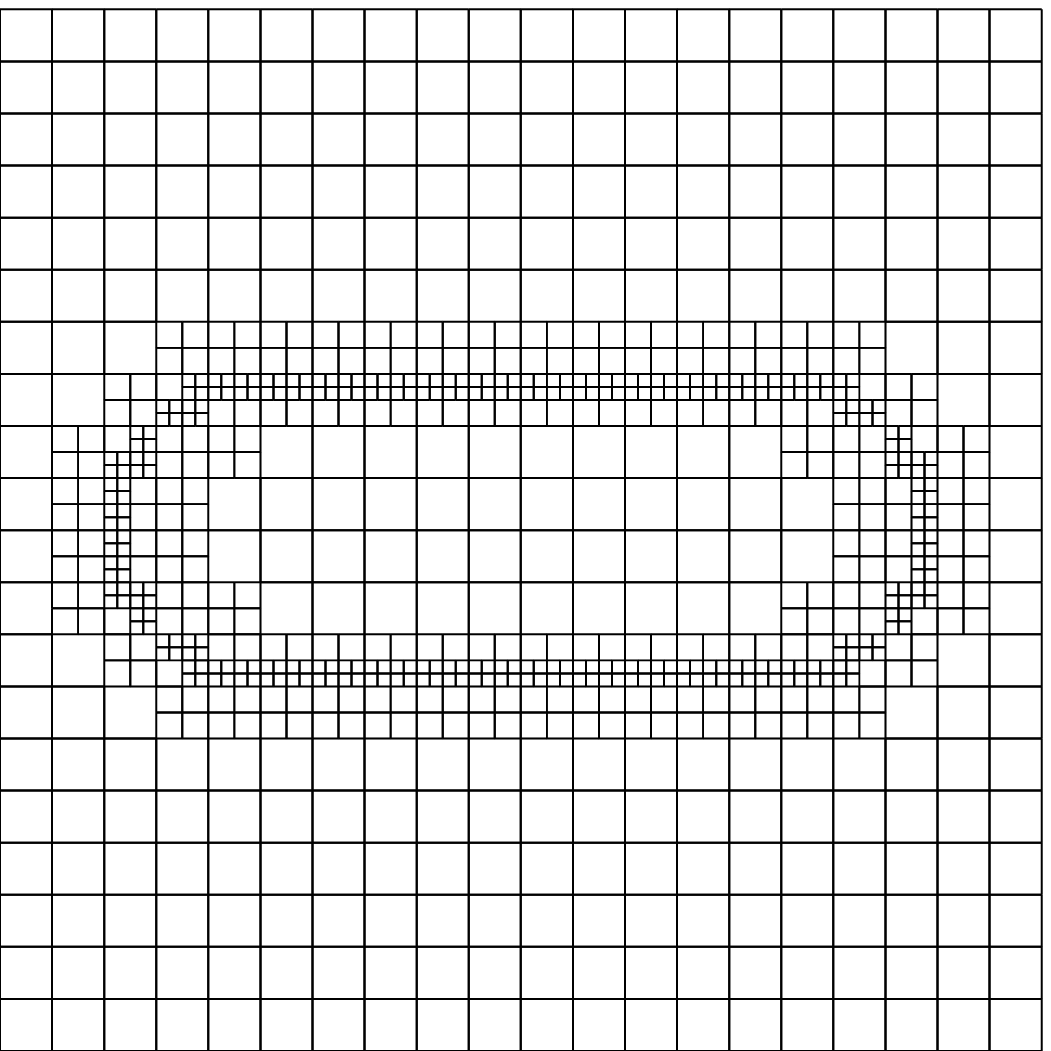}

}\subfloat[\label{fig:l2_geometry_projection}]{\includegraphics[width=0.25\textwidth]{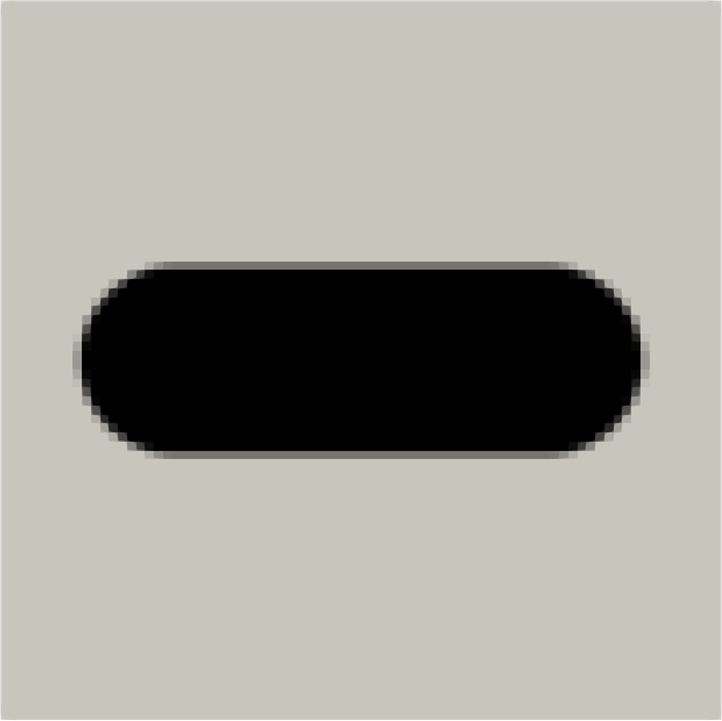}

}

\caption{Mesh refinement of a single component based on its effective density.
Top: first level of refinement; bottom: second level of refinement.
(a and e) effective density before refinement; (b and f) elements
marked for refinement; (c and g) refined mesh; and (d and h) effective
density after refinement.\label{fig:single-bar-refinement}}
\end{figure}

The AMR process for a structure made of a single solid bar (i.e.,
$\alpha=1$) is now illustrated for two levels of refinement. The
first level of refinement is shown in the top row of Fig.$\ $\ref{fig:single-bar-refinement}.
First, the composite density is computed on a coarse grid (Fig.\ \ref{fig:coarse_grid_projection}).
We use a sampling window radius $R$ corresponding to the circle that
circumscribes the element to compute its effective density $\hat{\rho}^{e}$
(in this 2-dimensional example, if the element $e$ size is $h^{e}$,
then $R=\sqrt{2}h^{e}/2$). We then mark elements for refinement based
on the refinement criterion of Eq.\ \ref{eq:eff-density-indicator},
shown in Fig.\ \ref{fig:coarse_grid_indicator}. Next, we refine
the marked elements by subdivision to produce the mesh shown in Fig.\ \ref{fig:refined_mesh_l1}.
We repeat the geometry projection on this refined mesh and re-compute
the element effective densities (Fig.\ \ref{fig:l1_grid_projection}).
With these new densities, we perform an additional level of refinement,
as shown in the bottom row of Fig.$\ $\ref{fig:single-bar-refinement}.
In this level, we not only mark elements for refinement based on the
refinement indicator, but also if they violate the single-level mesh-incompatibility
requirement. Additional levels of refinement follow the same procedure.
We note that, since the mesh refinement is based on the projected
density, the refined mesh provides a better resolution of the geometry
of the component, and the resolution of course increases with additional
levels of refinement, as can be seen by comparing Figs.$\ $\ref{fig:l1_grid_projection}
and \ref{fig:l2_geometry_projection}. 

\begin{figure}[h]
\subfloat[\label{fig:size-variable-one}]{\includegraphics[width=0.21\textwidth]{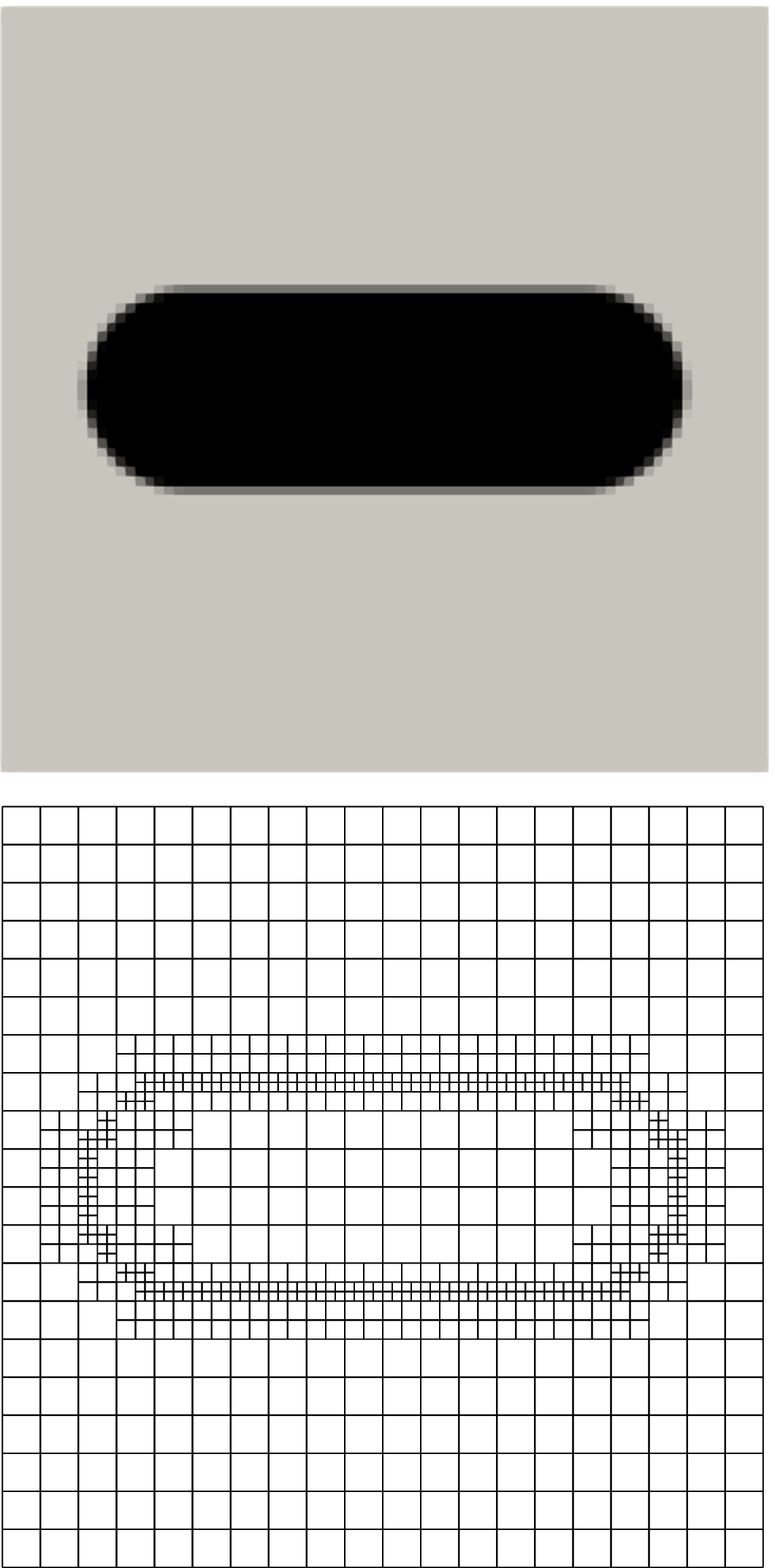}

}\hfill{}\subfloat[\label{fig:size-variable-m1}]{\includegraphics[width=0.21\textwidth]{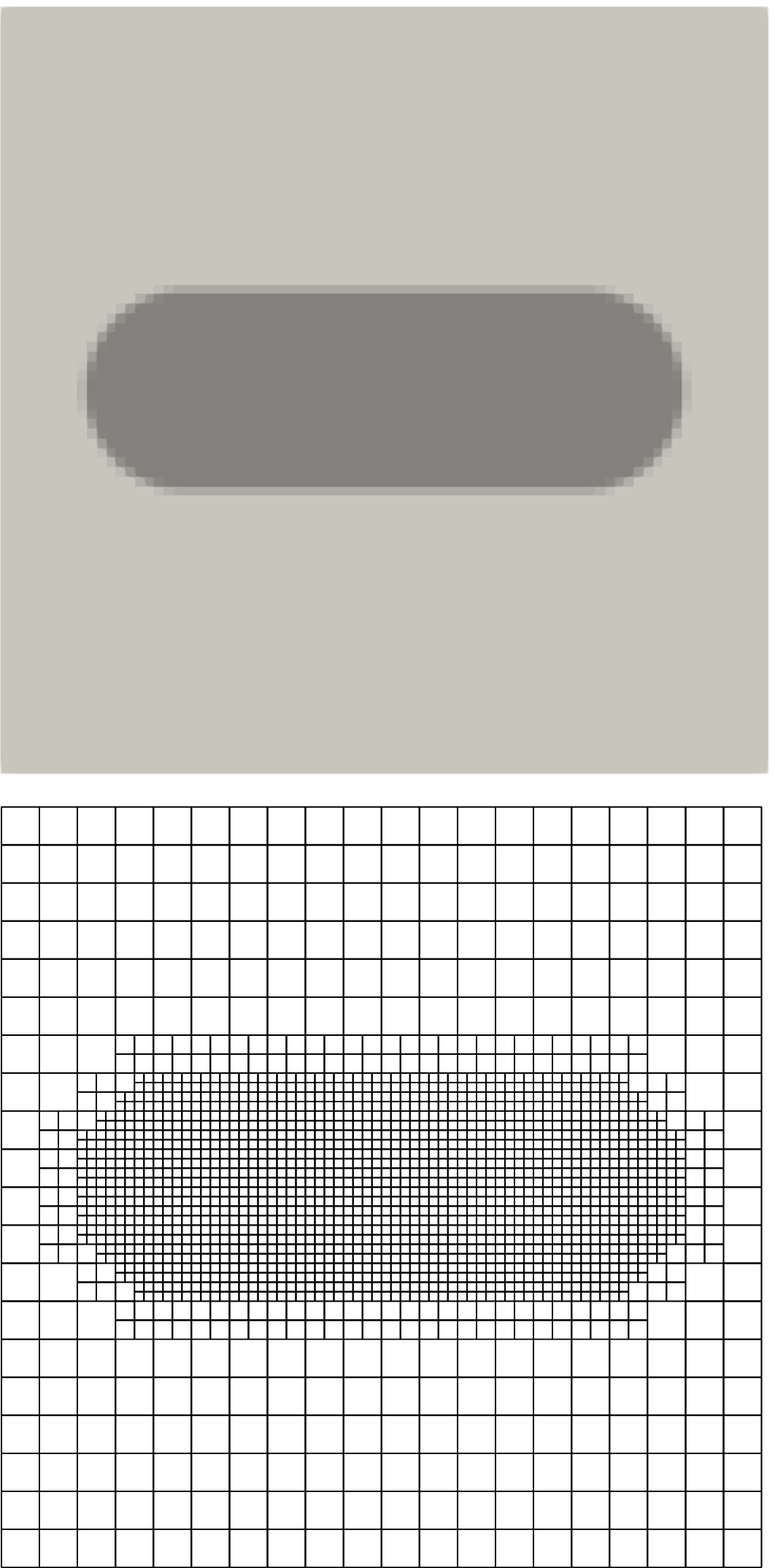}

}\hfill{}\subfloat[\label{fig:size-variable-m2}]{\includegraphics[width=0.21\textwidth]{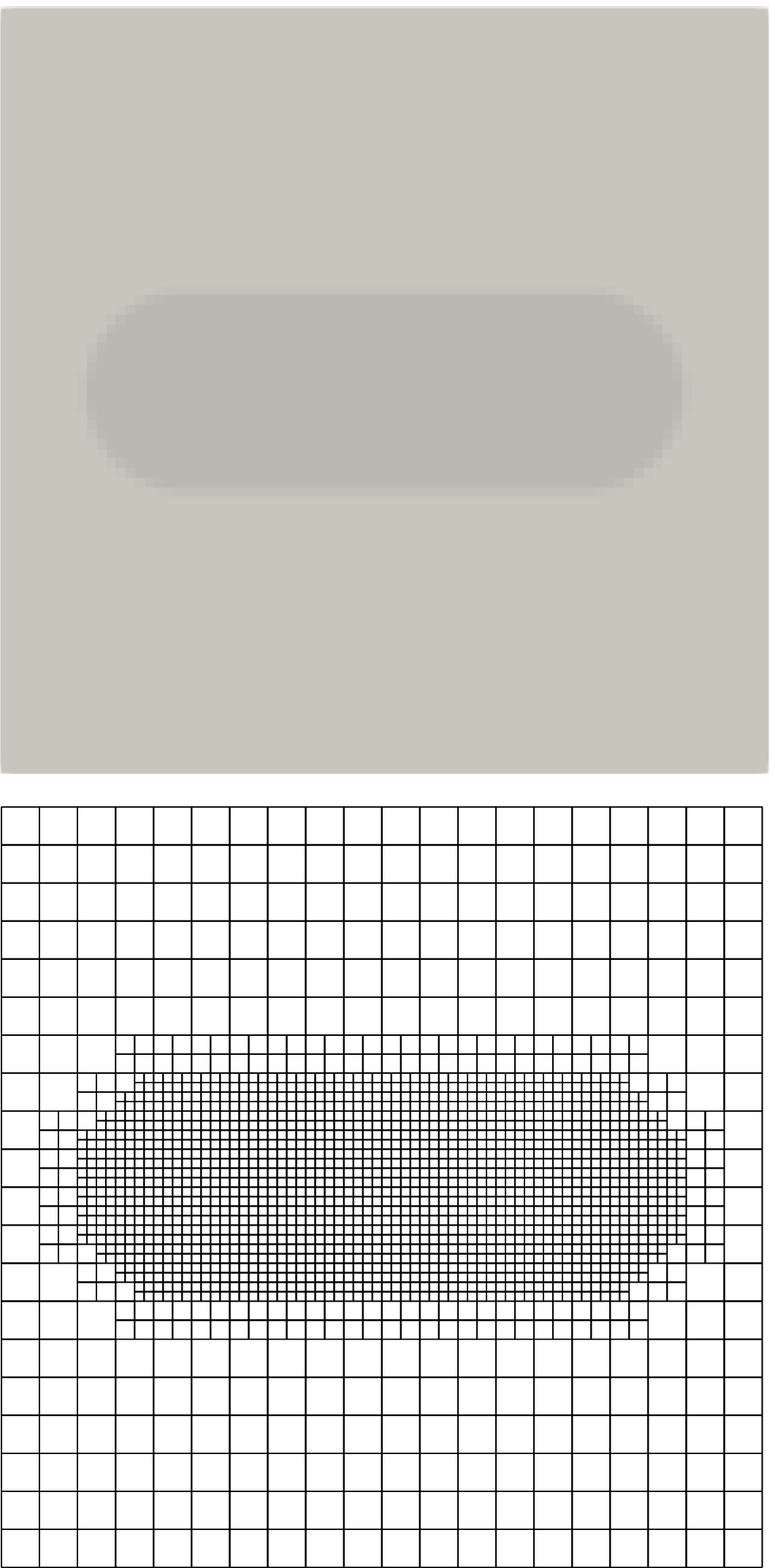}

}\hfill{}\subfloat[\label{fig:size-variable-zero}]{\includegraphics[width=0.21\textwidth]{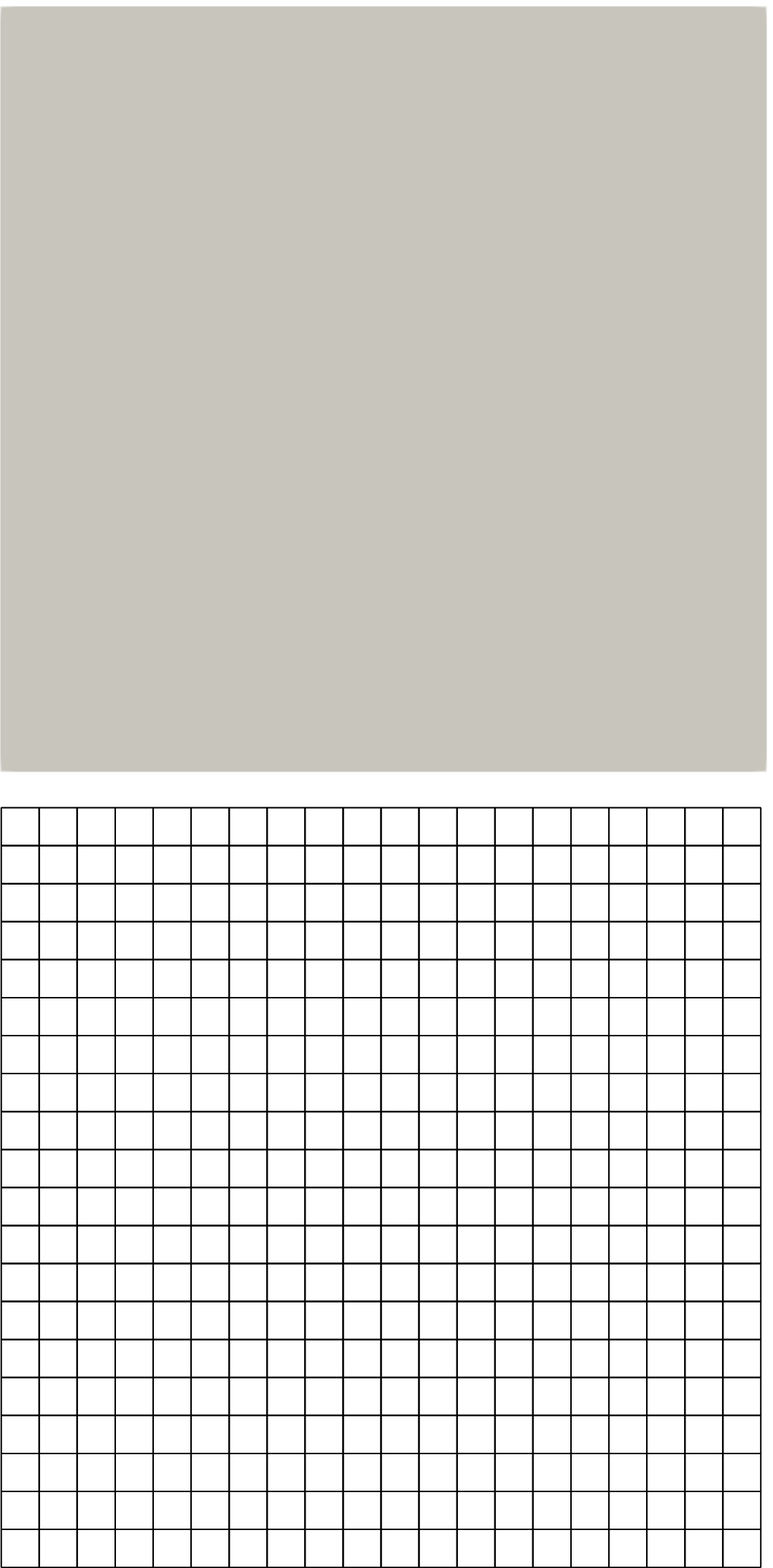}

}

\caption{Effect of the size variable on a two-level adaptive mesh refinement.
Effective density (top) and refined mesh (bottom) for bars with size
variables $\alpha$ of a) 1.0, b) 0.7, c) 0.4 and d) 0.0.\label{fig:size_variable_effect}}
\end{figure}

Since $\hat{\rho}^{e}$ is a function of the size variable of the
component, different values of the size variable $\alpha$ will result
in different refined meshes. Fig.\ \ref{fig:size_variable_effect}
demonstrates the effect of values of the size variable on the mesh
refinement. When $\alpha=1$, only elements near the boundary of the
bar are refined (Fig.\ \ref{fig:size-variable-one}). When $\alpha$
has an intermediate value (such that $\hat{\rho}^{e}<\rho_{th}$),
the entire region occupied by the bar is refined, which guarantees
elements with intermediate effective densities are always refined.
Finally, if $\alpha=0$, the bar is effectively removed from the design,
hence no mesh refinement is performed, as shown in Fig.\ \ref{fig:size-variable-zero}.

\begin{figure}[h]
\hfill{}\subfloat[\label{fig:nolayer}]{\begin{centering}
\includegraphics[width=0.4\textwidth]{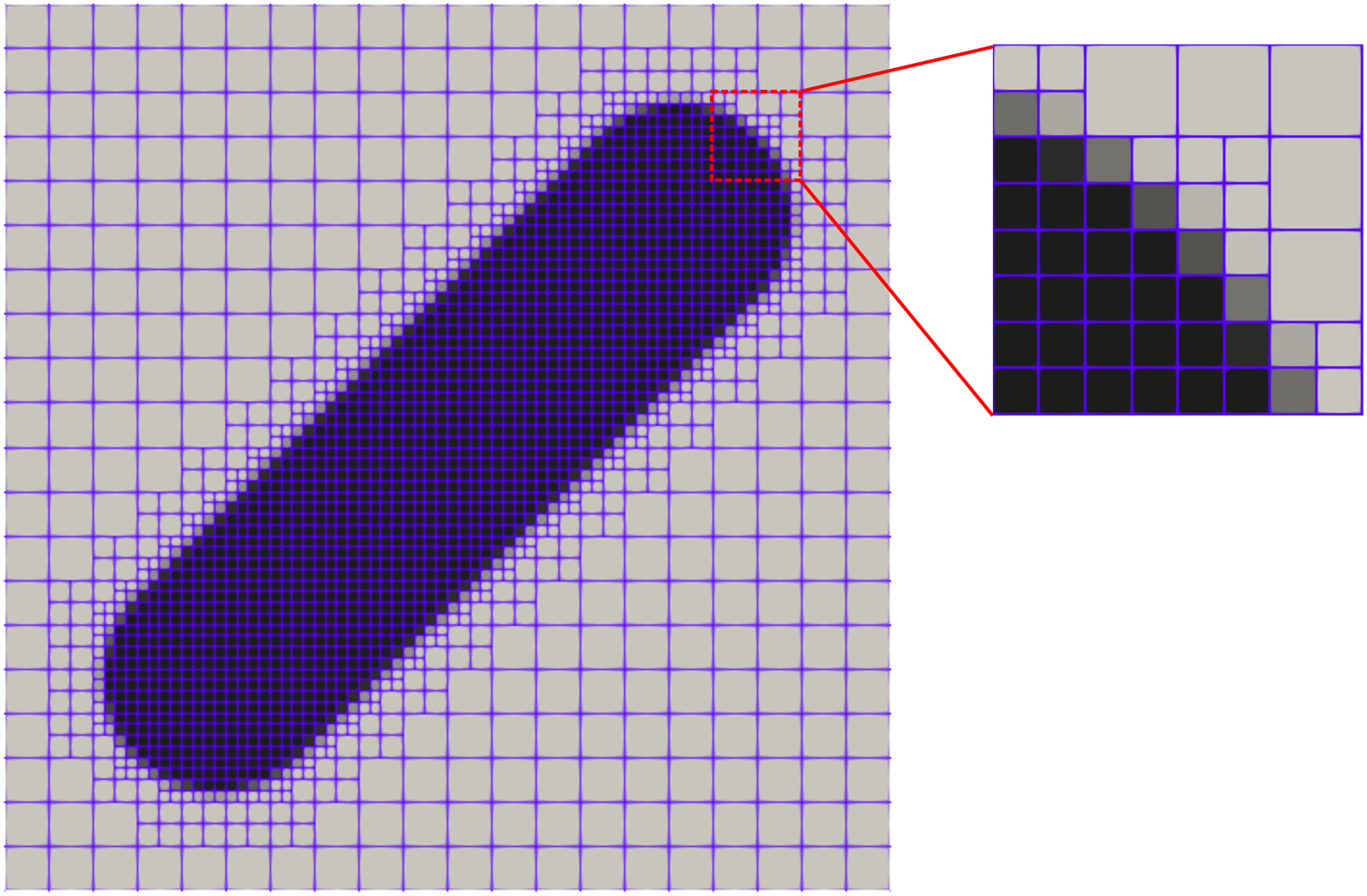}
\par\end{centering}
}\hfill{}\subfloat[\label{fig:withlayer}]{\begin{centering}
\includegraphics[width=0.4\textwidth]{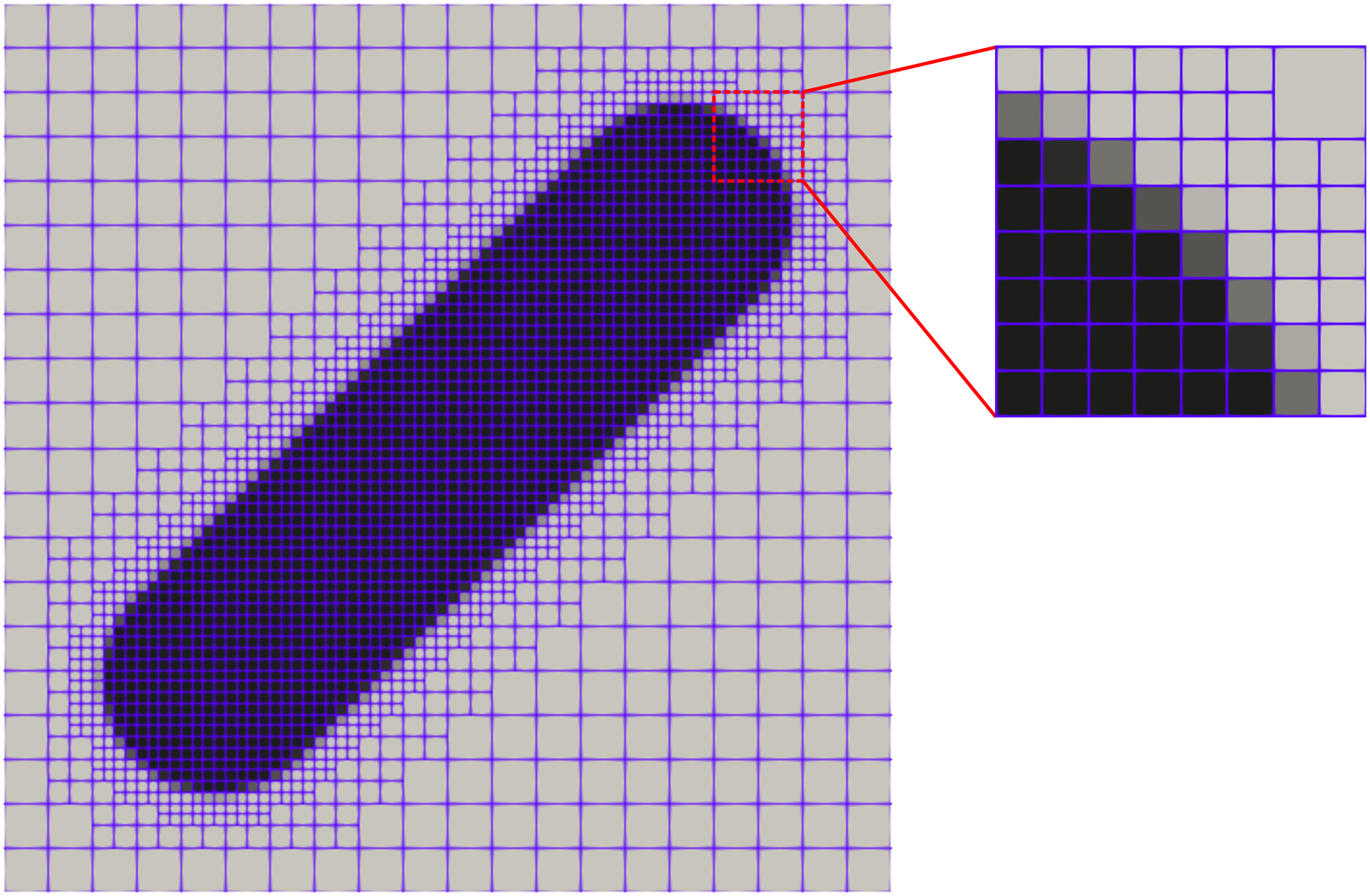}
\par\end{centering}
}\hfill{}

\caption{Element refinement along the component boundary: (a) Inconsistent
(b) Consistent.}
\end{figure}

\begin{figure}[h]
\center
\def\svgwidth{0.4\textwidth}     
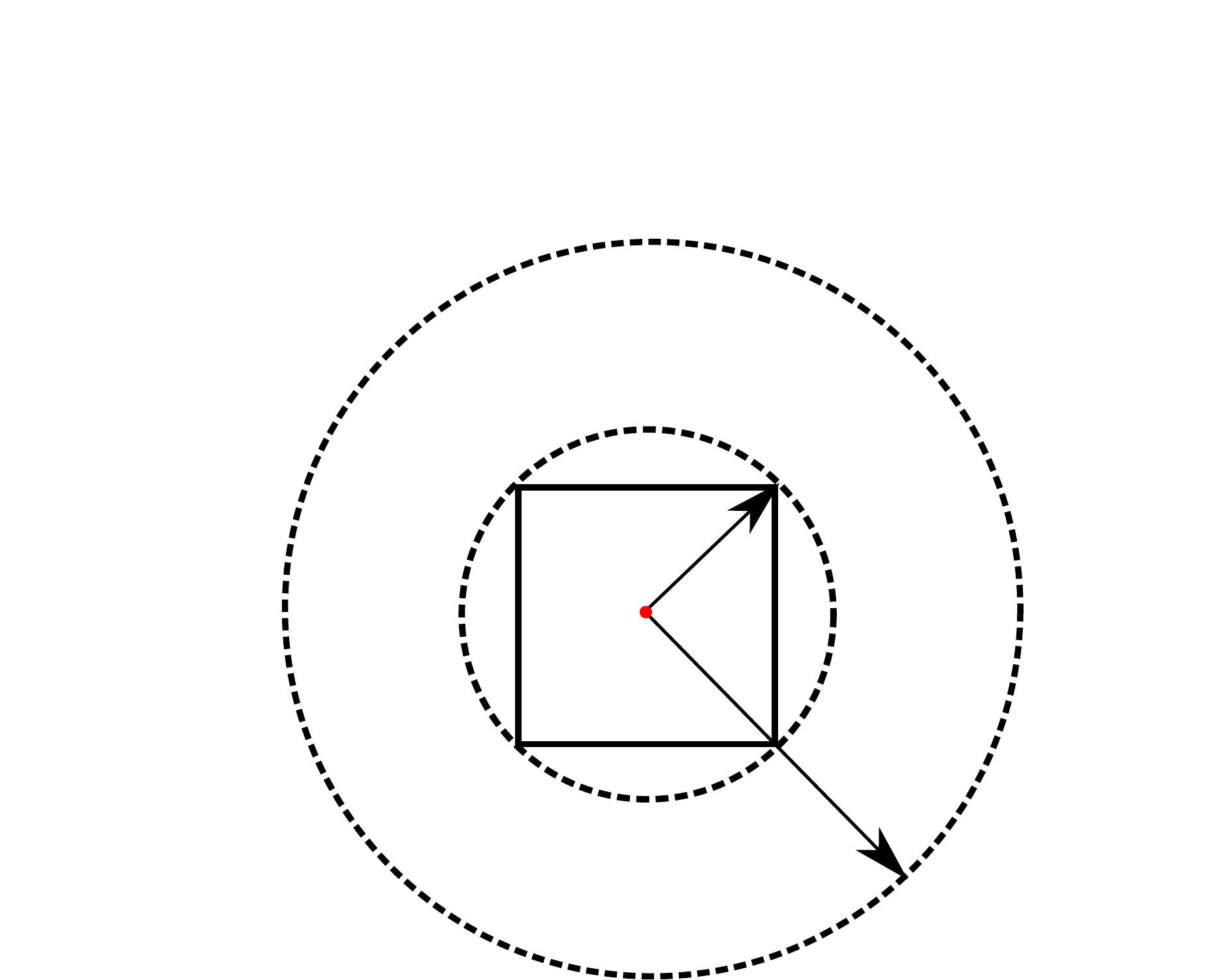

\caption{Consistent mesh size along boundary obtained by enlarging sampling
window.\label{fig:large_sampling_window}}

\end{figure}

The aforementioned refinement strategy may lead to an inconsistent
level of refinement (i.e., elements of different sizes) along the
boundary of the component. Fig.\ \ref{fig:nolayer} exemplifies this
situation, where we observe that the element size is not consistent
along the boundary of the bar due to the orientation of the bar. As
demonstrated in \cite{wang2010dynamic} for density-based topology
optimization, this inconsistency may lead to suboptimal designs. A
solution to this issue, proposed in \cite{wang2010dynamic}, is to
increase the size of the band of elements around the boundary that
are marked for refinement. It is straightforward to adapt this strategy
in our method by increasing the radius $R$ of the sampling window
used to compute the effective density. For example, Fig.\ \ref{fig:withlayer}
shows the refined mesh when we use a radius that equals twice the
radius of the circumscribing circle, i.e., $R=\sqrt{2}h^{e}$. In
effect, doubling $R$ is roughly equivalent to adding an additional
layer of elements to the band, as shown in Fig.$\ $\ref{fig:large_sampling_window}.
We therefore use this increased sample window radius for the mesh
refinement in all of our examples. We note, however, that when we
compute the effective density to obtain the ersatz material of Eq.\ \ref{eq:ersatz_material},
the circumscribing sampling window radius is still used to obtain
a more accurate geometry projection (note that, for a single bar,
the composite density equals the effective density of the bar).

\begin{figure}[h]
\begin{centering}
\subfloat[]{\def\svgwidth{0.25\textwidth}     
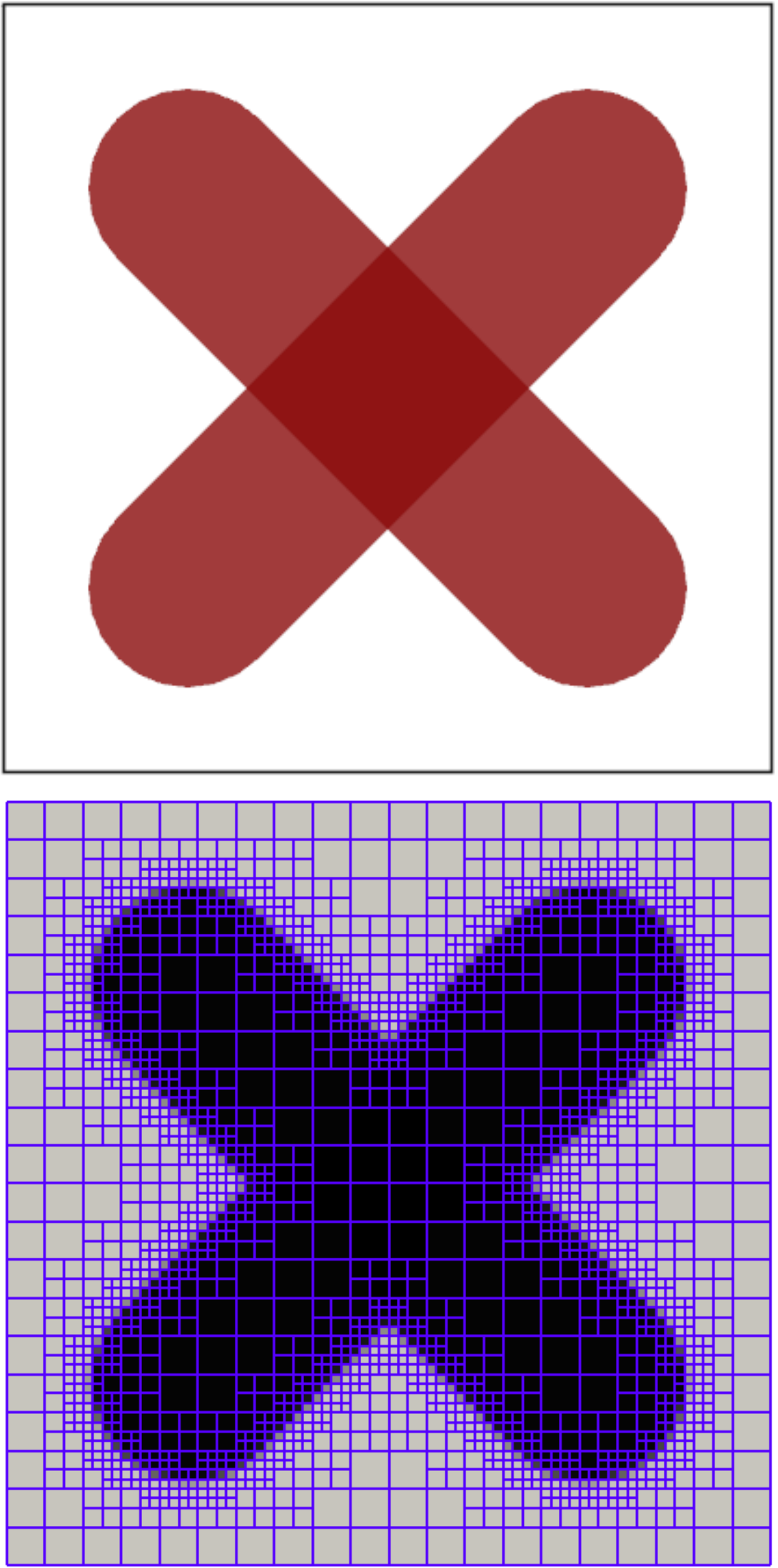

}\subfloat[]{\def\svgwidth{0.25\textwidth}     
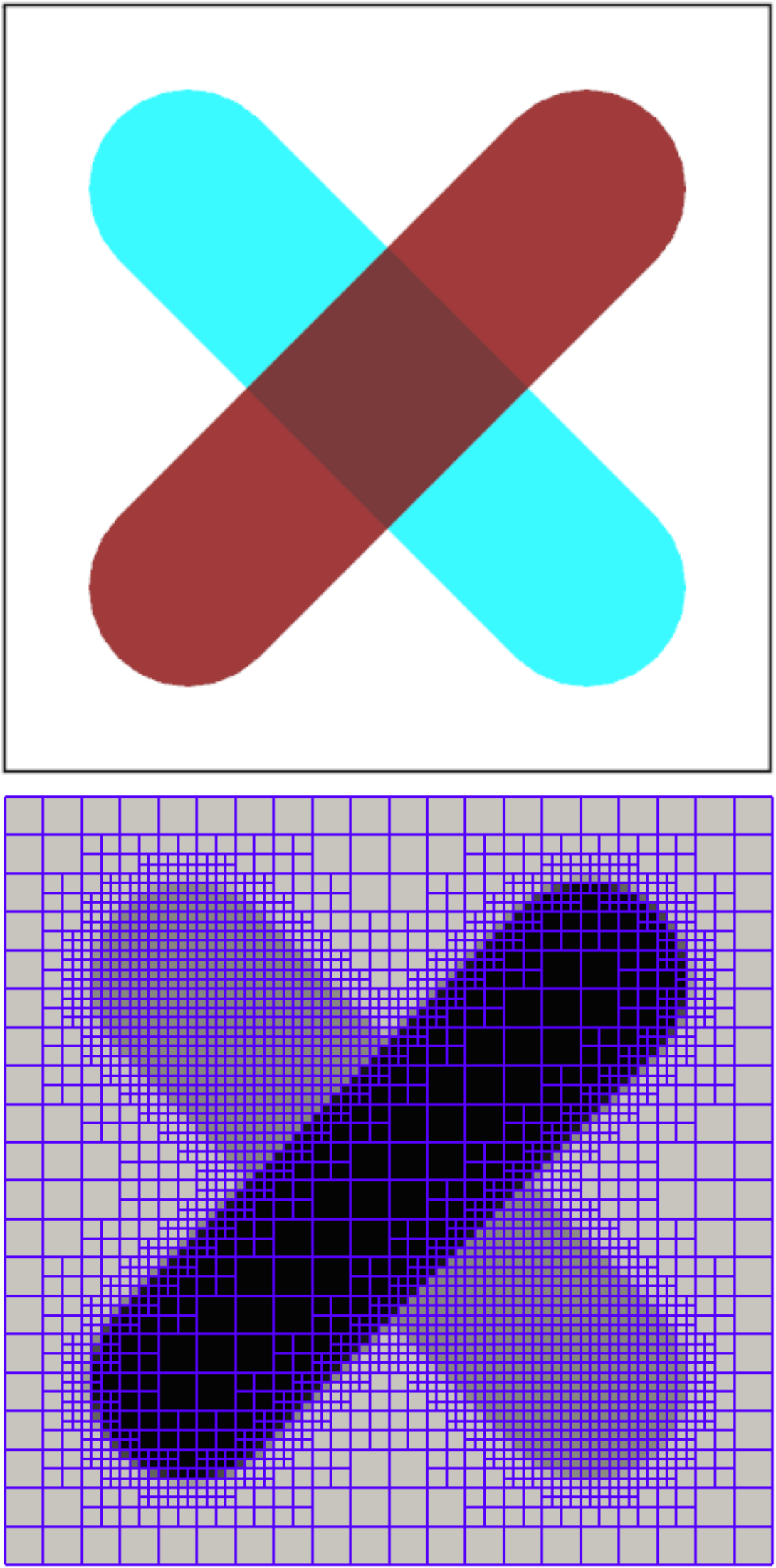

}\subfloat[]{\def\svgwidth{0.303\textwidth}     
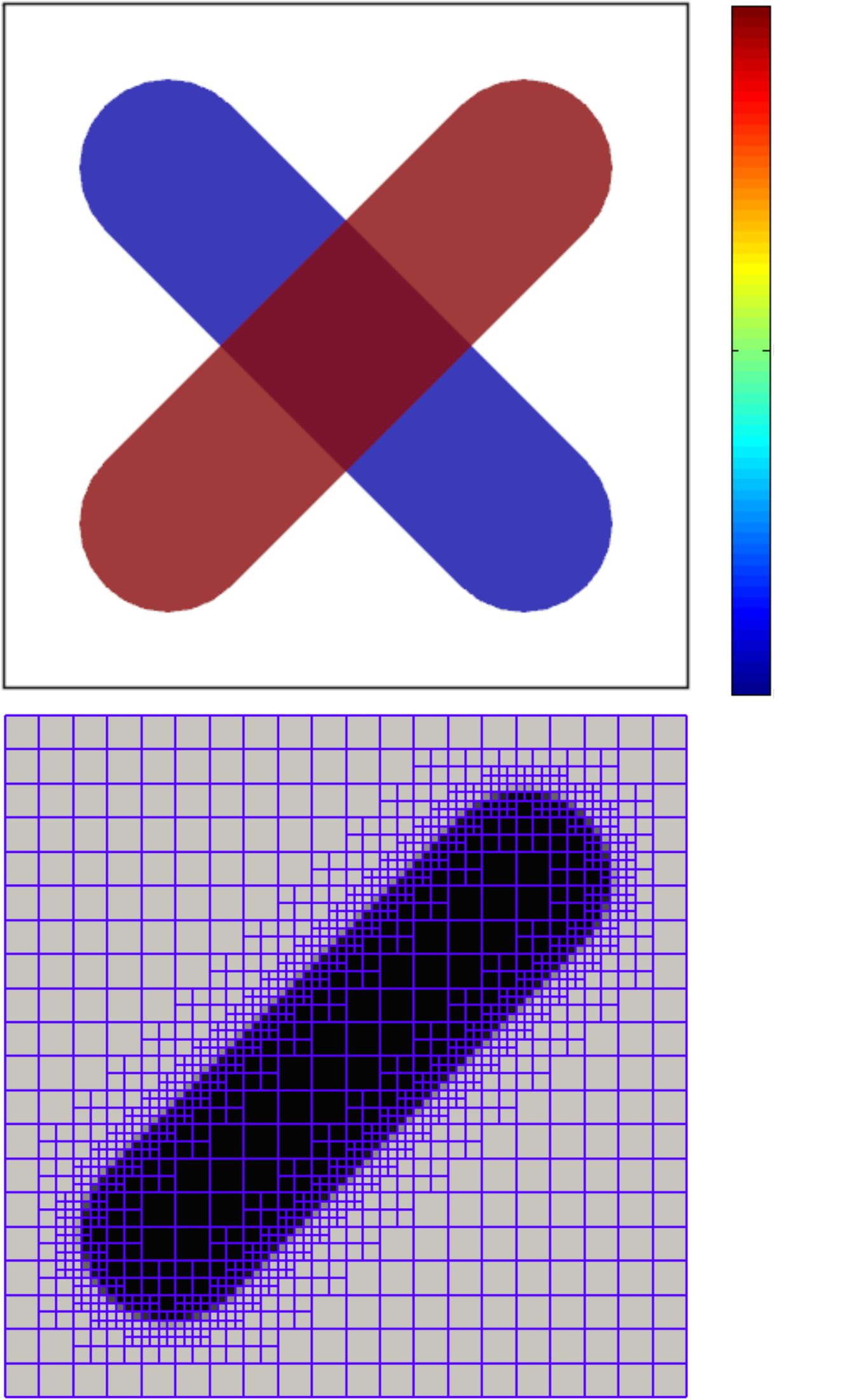

}
\par\end{centering}
\caption{Mesh refinement for two intersecting bars and the effect of the size
variable. Top: bars with solid colors indicating their penalized size
variable values. Bottom: refined mesh with two levels of refinement
and composite density. The size variable of the bar that runs from
top-left to right-bottom is (a) $\alpha=1$, (b) $\alpha=0.7$ and
(c) $\alpha=0$. \label{fig:two-bar-refinement}}

\end{figure}

The above discussion considered only a single component. When considering
multiple components, we simply impose the same refinement criterion
of Eq.$\ $\ref{eq:eff-density-indicator} on the composite density
of Eq.$\ $\ref{eq:composite-density-ks}, i.e.:
\begin{equation}
0<\tilde{\rho}_{ks}^{e}\leq\rho_{th}\lesssim1\label{eq:refinement-indicator-composite}
\end{equation}
When considering multiple bars, one important consideration is that
the lower-bound KS-function Eq.\ \ref{eq:LKS} underestimates the
true maximum. To see why this is important, we consider the case where
several bars overlap at the centroid of an element $e$, with one
of them having a size variable of $\alpha=1$ and the others having
$\alpha<1$. If we use the true maximum function to perform the Boolean
union, then the composite density would equal unity and therefore
the element would not be marked for refinement per Eq.$\ $\ref{eq:refinement-indicator-composite},
which is the behavior we desire because our main goal is to reduce
the mesh size as much as possible. However, when using the lower-bound
KS approximation of the maximum, then $\tilde{\rho}_{ks}^{e}<1$,
hence it is possible that $\tilde{\rho}_{ks}^{e}<\rho_{th}$ in spite
of the presence of the fully solid bar, and the element gets marked
for refinement. Therefore, we should not choose a threshold $\rho_{th}$
too close to unity to avoid this situation. The value $\rho_{th}=0.9$
is used in all the examples presented in this work and performed relatively
well. Fig.$\ $\ref{fig:two-bar-refinement} illustrates the use of
the composite density as a refinement criterion for multiple bars,
with several values of the size variable of one of the bars. As seen
in the figure, when both bars have a size variable of unity, the entire
region inside of their union has a coarse mesh; when one of the bars
has an intermediate value, the portion of that bar that does not intersect
the other bar results in a fine mesh; and when the size variable of
that bar is zero, it has no effect on the mesh refinement, because
in that case the composite density equals the effective density of
the other bar.

\section{Optimization problems and computer implementation\label{sec:Optimization-problems-and}}

To demonstrate the effectiveness of the proposed AMR method with discrete
geometric components, we consider two optimization problems: 1) a
minimum compliance problem; and 2) a stress-constrained problem. For
the compliance-based problem, we consider the minimization of the
structural compliance subject to a volume constraint, formulated as:

\begin{align}
\text{\ensuremath{\min_{\mathbf{z}}\;}C(\ensuremath{\mathbf{u}})} & :=\int_{s^{t}}\mathbf{u}\cdot\mathbf{t}\ da\label{eq:objective_function}\\
\text{subject to:}\nonumber \\
v_{f}(\mathbf{z}) & \leq v_{f}^{\star}\label{eq:constraints}\\
\mathsf{a}(\mathbf{u},\mathbf{v}) & =\mathsf{l}(\mathbf{v})\quad\mathbf{u},\mathbf{\forall v}\in\mathcal{U}_{ad}\label{eq:bi-linear}\\
\text{\ensuremath{\underbar{\ensuremath{\mathbf{z}}}}}\leq\mathbf{z} & \leq\bar{\mathbf{z}}\label{eq:dv_bounds}
\end{align}
where C denotes the structural compliance and the volume fraction
$v_{f}$ of the structure is computed via

\begin{equation}
v_{f}(\mathbf{z}):=\frac{\mid\hat{\omega}\mid}{\mid\Omega\mid}=\frac{1}{\mid\Omega\mid}\int_{\Omega}\tilde{\rho}_{ks}(\mathbf{x},\mathbf{z},0,R)\thinspace dv\label{eq:volume-fraction}
\end{equation}
In this equation, $\hat{\omega}=\bigcup_{i}\omega_{i}\subset\Omega$
denotes the entire structure and $\Omega$ denotes the design envelope.
The constraint limit on the volume fraction is denoted by $v_{f}^{\star}$.
In Eq.\ \ref{eq:bi-linear}, $\mathsf{a}$ is the energy bilinear-form
and $\mathsf{l}$ is the load linear-form given by

\begin{align}
\mathsf{a}\left(\mathbf{u},\mathbf{v}\right) & :=\int_{\Omega}\nabla\mathbf{v}\cdot\mathbb{C}(\mathbf{x},\mathbf{z},s,R)\nabla\textbf{\ensuremath{\mathbf{u}}}\thinspace dv\label{eq:energy-bilinear}\\
\mathsf{l}\left(\textbf{\ensuremath{\mathbf{v}}}\right) & :=\int_{S^{t}}\textbf{\ensuremath{\mathbf{v}}}\cdot\mathbf{t}\ da\label{eq:load-linear}
\end{align}
where $\mathbf{u},\mathbf{v}\in\mathcal{U}_{ad}$ are the displacement
and test functions respectively, and $\mathcal{U}_{ad}:=\{\mathbf{u}\in H^{1}(\Omega):\mathbf{u}=\mathbf{0}\text{ on }S^{u}\}$
is the set of admissible displacements. The boundary of the structure
$\partial\hat{\omega}$ is composed of the zero traction boundary
$S^{0}$, the displacement boundary $S^{u}$ and the nonzero traction
boundary $S^{t}$. We assume the boundaries $S^{0}$, $S^{u}$ and
$S^{t}$ are not empty, and that $S^{u},S^{t}\subset\partial\Omega$
are design-independent. The lower and upper bounds on the design variables
$\mathbf{z}$ are denoted by $\text{\ensuremath{\underbar{\ensuremath{\mathbf{z}}}}}$
and $\bar{\mathbf{z}}$ respectively. To improve the convergence of
the optimization, all design variables are scaled to the interval
$[0,1]$. Therefore, Eq.\ \ref{eq:dv_bounds} is modified as:

\begin{equation}
\text{\ensuremath{\underbar{\ensuremath{\mathbf{b}}}}}\leq\hat{\mathbf{z}}\leq\bar{\mathbf{b}}\label{eq:scaled-dv-bounds}
\end{equation}
where $\text{\ensuremath{\underbar{\ensuremath{\mathbf{b}}}}}$, and
$\bar{\mathbf{b}}$ are the scaled bounds on the scaled design variables
$\hat{\mathbf{z}}$. The scaling strategies for the optimization of
bar, plate, curved-plate and supershape components are detailed in
\cite{norato2015geometry}, \cite{zhang2016geometry}, \cite{zhanggeometry}
and \cite{norato2018topology} respectively. In Eqs.\ \ref{eq:volume-fraction}
and \ref{eq:energy-bilinear}, $R$ is set to $\sqrt{2}h^{e}/2$ for
the 2-dimensional examples and $\sqrt{3}h^{e}/2$ for the 3-dimensional
examples, corresponding to the radii of the circle and the sphere
that circumscribe the square and cube elements respectively. 

For the stress-based optimization problem, we consider volume minimization
subject to a constraint on the maximum stress of the structure:

\begin{align}
\text{\ensuremath{\min_{\mathbf{z}}\ }} & v_{f}(\mathbf{z})\label{eq:vol_obj}\\
\text{subject to:}\nonumber \\
\max_{\mathbf{x}\in\hat{\omega}}\,\sigma(\mathbf{x},\mathbf{z}) & \leq\sigma^{*}\label{eq:stress_cons}\\
\mathsf{a}(\mathbf{u},\mathbf{v}) & =\mathsf{l}(\mathbf{v})\quad\mathbf{u},\mathbf{\forall v}\in\mathcal{U}_{ad}\label{eq:bi-linear-2}\\
\text{\ensuremath{\underbar{\ensuremath{\mathbf{b}}}}}\leq\hat{\mathbf{z}} & \leq\bar{\mathbf{b}}\label{eq:dv_bounds-2}
\end{align}
where $\sigma(\mathbf{x},\mathbf{z})$ is the von Mises stress at
$\mathbf{x}$ and $\sigma^{*}$ is the maximum allowable stress. For
the sake of brevity, we refer readers to \cite{zhang2017stress}
for details on the stress constraint of Eq.\ \ref{eq:stress_cons}
and on the optimization strategies employed to solve the stress-based
topology optimization problem with geometric components.

We use the Method of Moving Asymptotes (MMA) in \cite{svanberg1987method,svanberg1995globally,svanberg2002class}
to solve the foregoing optimization problems. In specific, Caterpillar's
implementation of MMA based on \cite{svanberg1995globally} is used,
where the subproblems are solved using the L-BFGS-B method of \cite{byrd1995limited}
and the coefficient $M$ is set to 10 unless otherwise stated. In
addition, we impose a uniform move limit on all scaled design variables
as in \cite{norato2015geometry} to improve the convergence of the
optimization. For a scaled design variable $\hat{z}$ at iteration
$I$, a move limit $0<m\leq1$ is imposed as

\begin{equation}
\max\left(\text{\ensuremath{\underbar{b}}},\hat{z}^{I-1}-m\right)\leq\hat{z}^{I}\leq\min\left(\bar{b},\hat{z}^{I-1}+m\right)\label{eq:move-limit}
\end{equation}

Our code is implemented using the \texttt{deal.II} finite element
library \cite{BangerthHartmannKanschat2007,dealII90}. The library
provides mesh refinement utilities, including the handling of the
resulting hanging nodes in the analysis, greatly facilitating the
implementation of our proposed AMR scheme. Moreover, to ensure the
efficiency of our method, we parallelize the AMR, geometry projection,
element assembly, solution of the linear system and sensitivity calculation
using the parallel data structures and linear algebra software provided
with the library.

\begin{figure}[h]
\begin{centering}
\includegraphics[width=0.6\textwidth]{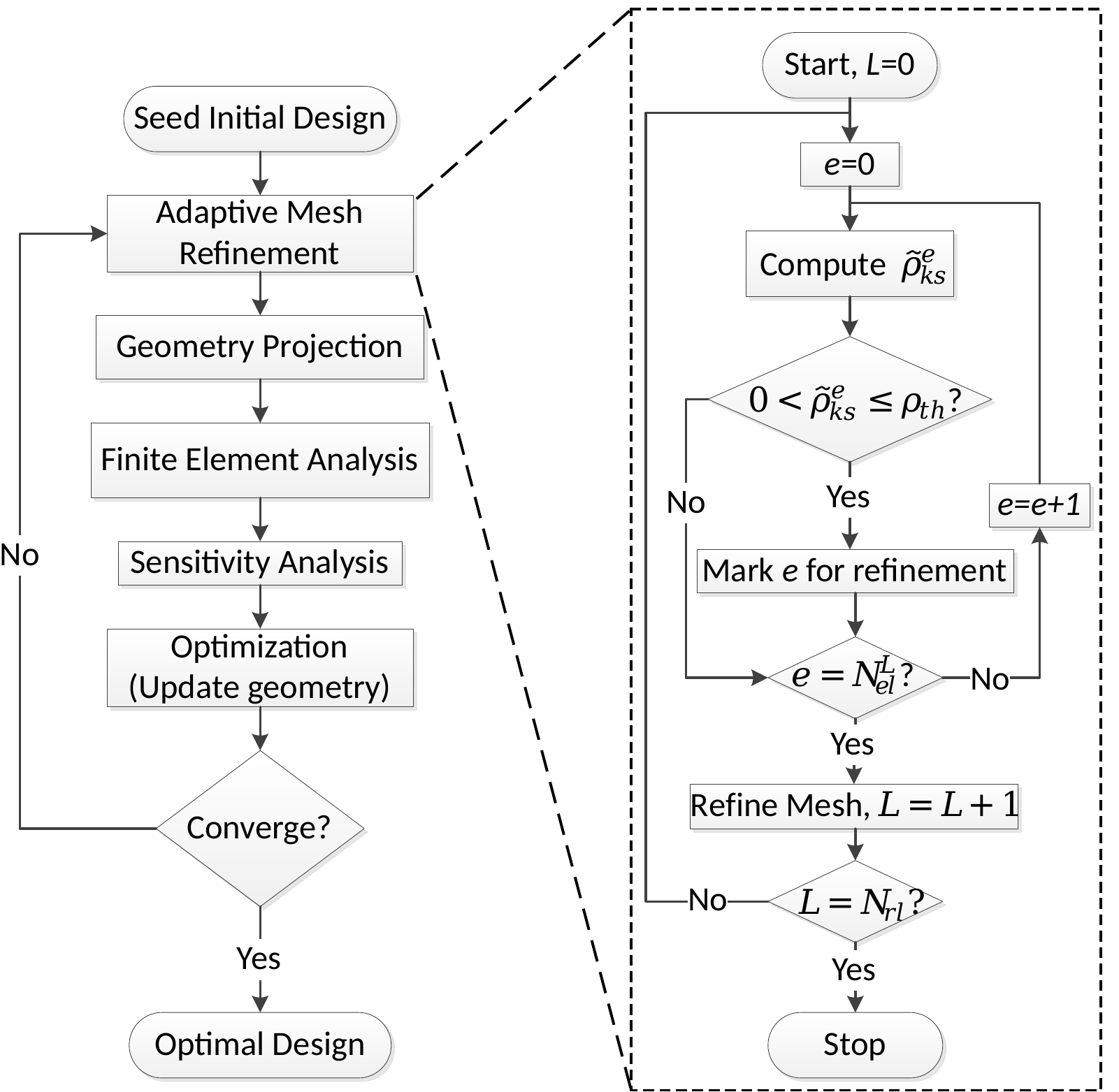}\caption{Flowchart of the optimization and the AMR.\label{fig:flow-chart}}
\par\end{centering}
\end{figure}

The topology optimization with AMR is summarized in the flowchart
of Fig.\ \ref{fig:flow-chart}. Within each outer optimization loop,
the mesh refinement always starts from the coarsest mesh level $L=0$.
As illustrated in the previous Section, the next mesh refinement level
is obtained by marking elements for refinement based on the refinement
criterion and the single-level mesh-incompatibility requirement; marked
elements are then refined by subdivision, and the process repeated
until the desired refinement level $N_{rl}$ is attained. Since we
start from a coarse mesh at every outer optimization loop, when the
design changes it is entirely possible that regions that were refined
in previous design iterations are no longer refined in the current
one, and thus have a coarse mesh. 

We note that in density-based methods it is necessary to re-define
the set of design variables every time the mesh is adapted, since
the design representation is tied to the mesh; moreover, the optimization
functions corresponding to different mesh discretizations are mapped
from previous optimization iterations to the current to avoid having
to restart MMA (cf., \cite{Troya2018}). In the topology optimization
with geometric components that we consider in this work, these strategies
are unnecessary because the representation of the design is independent
of the mesh. We also note that the geometry projection-based indicator
of Eq.$\ $\ref{eq:refinement-indicator-composite} provides a straightforward
way of performing the AMR at each iteration of the optimization, which,
as demonstrated in \cite{wang2010dynamic} and mentioned in Section
\ref{sec:Introduction}, leads to better designs.

The computational bottleneck of the topology optimization is the solution
of the system of linear equations. For modest size problems, direct
solvers can be used to efficiently solve this linear system, and they
have the advantage of being less sensitive to the large condition
number arising by the high-contrast between void and solid regions
in topology optimization. However, for the large-scale 3-dimensional
problems presented later in this paper, the use of direct solvers
becomes prohibitive because of their computational and storage requirements.
Iterative solvers such as the Conjugate Gradient (CG) method, on the
other hand, have relatively lower computational cost and storage requirements
for large problems. To achieve fast convergence for iterative solvers,
a good preconditioner is necessary. Unfortunately, popular preconditioners
such as incomplete Cholesky and Jacobi are very sensitive to the condition
number of the system. This problem is exacerbated in the topology
optimization with geometric components, because when the components
do not form a connected load path between the points of application
of the loads and the displacement boundary conditions, the conditioning
of the system substantially worsens. Moreover the conditioning of
the system also worsens with an increased number of elements. As the
condition number increases, the solution will require more iterations
to converge, which makes the solution time infeasible. Therefore,
it is desired to have a preconditioner for the iterative solution
for large-scale problems that is less sensitive to the condition number
and that scales well with the problem size. As shown in \cite{amir2014multigrid},
the Geometric MultiGrid (GMG) preconditioner with CG provides superior
performance for the solution of the linear systems arising in topology
optimization. A framework of the GMG preconditioner on adaptively
refined meshes is proposed in \cite{janssen2011adaptive} and implemented
in the \texttt{deal.II} library; we employ it for the 3-dimensional
examples to speed up the solution process.

\section{Examples}

We present several numerical examples to demonstrate the proposed
method. For all the examples, we consider all geometric components
to be made of a homogeneous, isotropic, linear elastic material with
Young\textquoteright s modulus $E=1\text{E}5$ and Poisson\textquoteright s
ratio $\nu=0.3$. We employ bar and plate components for the 2-dimensional
and 3-dimensional design problem respectively. The design envelope
is discretized with bilinear quadrilateral elements for the 2-dimensional
problems, and with trilinear hexahedral elements for the 3-dimensional
problems. As aforementioned, we employ a direct solver for the 2-dimensional
problems, and an iterative solver (preconditioned conjugate gradient)
for the 3-dimensional problems. We impose the lower bound $\rho_{\min}=1\text{E-}4$
on the void region (cf.$\ $Eq.$\ $\ref{eq:ersatz_material}). The
optimization is stopped when the absolute value of the maximum change
in the scaled design variables between the current iteration and two
previous consecutive iterations is less than $5\text{E-}3$. All the
2-dimensional examples are performed using a single thread on an 8-core
3.60GHz Intel Core i7-7820X processor. The machine information for
the 3-dimensional problems will be described separately for each example.

\subsection{2-dimensional MBB beam design for compliance minimization}

\begin{figure}[h]
\center
\def\svgwidth{0.6\textwidth}     
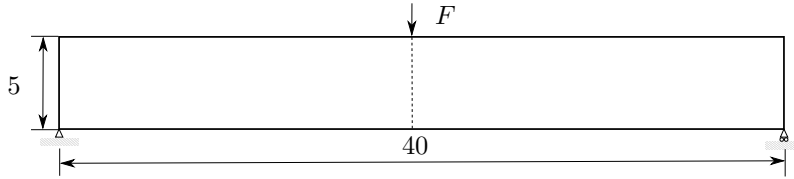\caption{MBB beam design problem definition\label{fig:MBB-beam-dims}}

\end{figure}

\begin{figure}[h]
\center
\def\svgwidth{0.6\textwidth}     
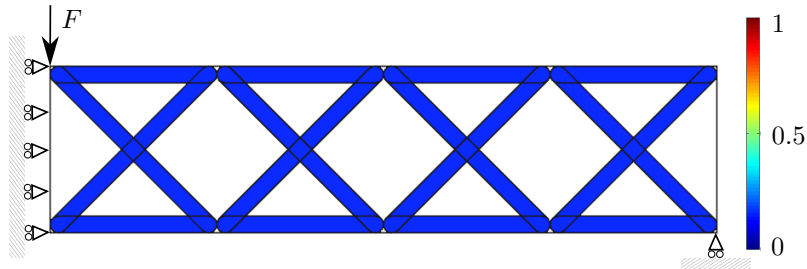

\caption{Initial design, geometry, loads and boundary conditions for the half-MBB
beam design problem. Color denotes the penalized size variable $\alpha^{s}$.\label{fig:MBB-initial}}
\end{figure}

We consider the well-known Messerchmitt-B\"{o}lkow-Blohm (MBB) beam design
problem shown in Fig.\ \ref{fig:MBB-beam-dims} for the first example.
The optimization problem is the minimization of the structural compliance
subject to a volume fraction constraint $v^{*}=0.3$ (Eq.\ \ref{eq:volume-fraction}).
Since the loading and boundary conditions are symmetric with respect
to the center plane of the beam (indicated by a dashed line in Fig.\ \ref{fig:MBB-beam-dims}),
we only model the right-hand side of the beam. The design envelope,
initial design, external loading and boundary conditions for the half-beam
are shown in Fig.\ \ref{fig:MBB-initial}. The applied load $F$
has a magnitude of 10. The move limit $m$ in Eq.\ \ref{eq:move-limit}
is set to 0.05. In this example, the coarsest mesh level, corresponding
to $L=0$, is a uniform mesh with element size $h_{c}=0.25$. The
target refinement level $N_{rl}$ is set to 2 so that the finest mesh
element has size $h=0.0625$ after the AMR. The strategy described
in Section \ref{sec:Adaptive-mesh-refinement} and Fig.\ \ref{fig:flow-chart}
is used to generate the adaptively refined mesh at each optimization
iteration. 

\begin{figure}[h]
\begin{centering}
\subfloat[\label{fig:mbb-AMR-iteration}]{\includegraphics[width=0.475\textwidth]{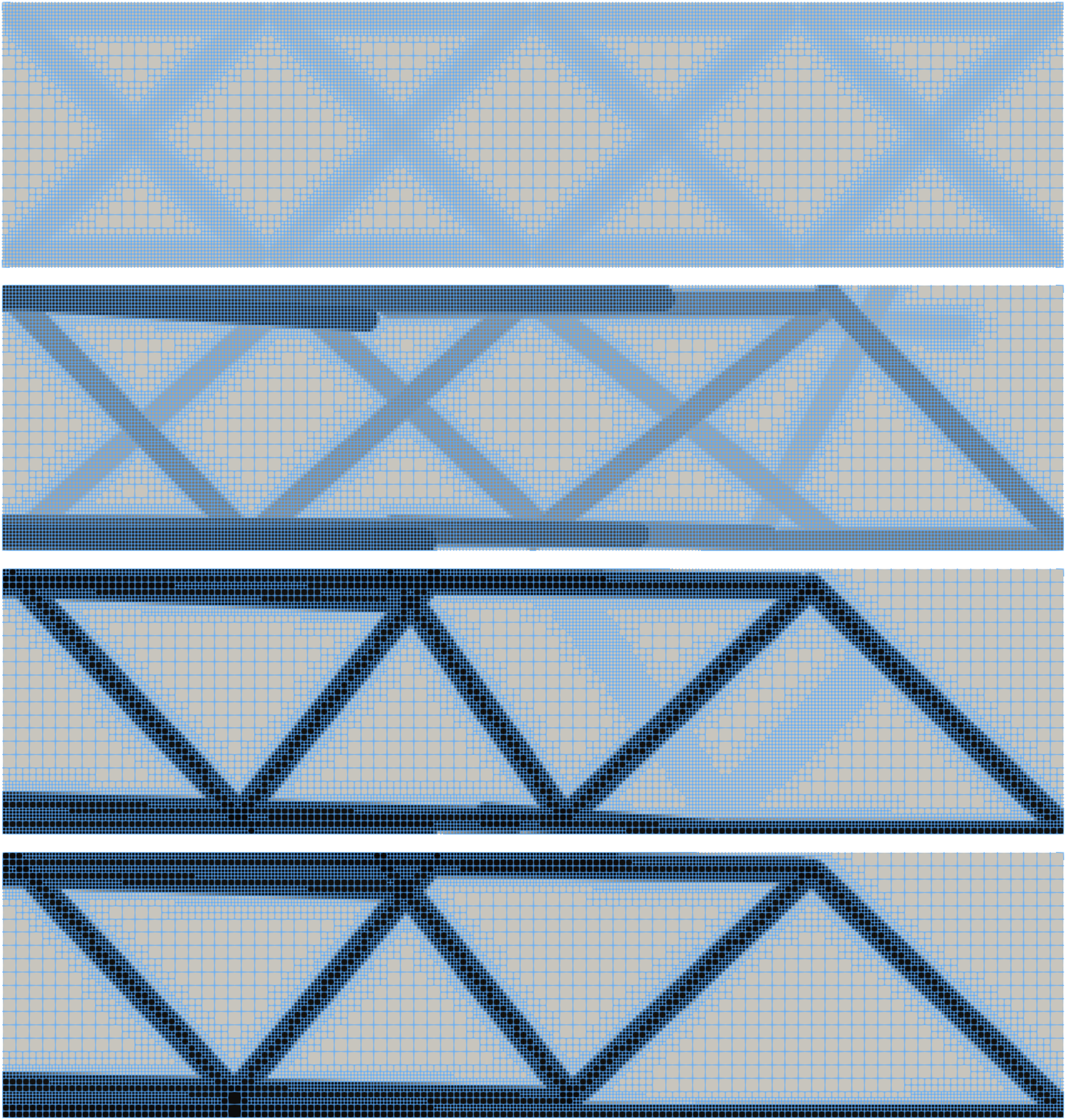}

}\subfloat[\label{fig:mbb-full-iteration}]{\includegraphics[width=0.475\textwidth]{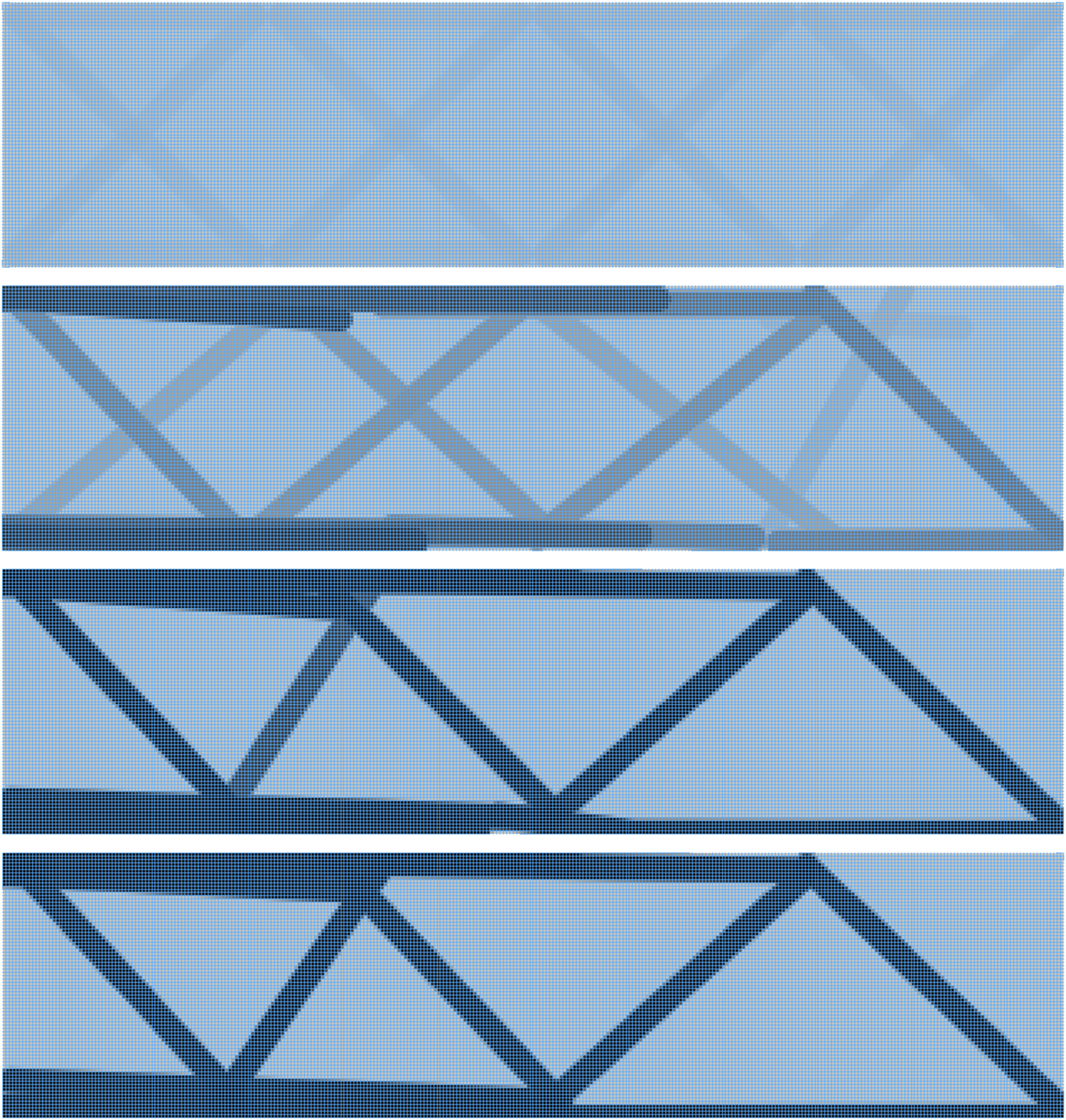}

}
\par\end{centering}
\caption{MBB beam design optimization with AMR (left) and full-resolution (right)
meshes. From top to bottom: initial design, iteration 10, iteration
30, and last iteration (50 for AMR and 42 for full-resolution mesh).
Compliance values for optimal designs are $C=0.534$ for AMR and $C=0.539$
for full-resolution mesh.}
\end{figure}

The AMR mesh and composite density for the initial design, iterations
10 and 30, and the last iteration (50) are plotted in Fig.\ \ref{fig:mbb-AMR-iteration}.
The design iterates show that the mesh is adaptively refined as the
optimization progresses. The region with intermediate density is meshed
with a fine mesh, while regions with zero or near-unity density are
meshed using coarse elements. The optimization converges to a design
with compliance $C=0.534$ in 166 seconds and the mesh in the last
iteration consists of 10,756 elements. We perform the optimization
for the same problem on a uniform fine mesh with element size $h=0.0625$
and a total 25,600 elements. This full-resolution mesh is obtained
by globally refining the coarsest mesh $N_{rl}$ times before the
optimization (we use this same procedure for the full-resolution meshes
in this and all of the following examples). The mesh and composite
density for the initial design, iterations 10 and 30, and the last
iteration (42) are plotted in Fig.\ \ref{fig:mbb-full-iteration}.
The entire optimization takes 171 seconds to converge to a design
with $C=0.539$, which is fairly similar to the one obtained with
AMR. Moreover, the two designs are fairly similar. The average solution
time per iteration is 3.3 seconds for AMR and 4.1 seconds with the
full-resolution mesh. The timing information for the last optimization
iteration is shown in Tab.\ \ref{tab:time}. 

In this example, using AMR does not significantly reduce the the total
optimization time, because the mesh refinement itself takes about
40\% of the total time for this small problem. Nevertheless, AMR successfully
halves the number of mesh elements, which reduces the computational
burden for any calculation that is a function of the number of mesh
elements, including the geometry projection and its sensitivities,
the global stiffness matrix computation and the finite element analysis.
As a result, using AMR still outperforms using the full-resolution
mesh despite the overhead time associated with the AMR. As we demonstrate
in other examples, however, the savings achieved by using AMR is significant
for large-scale problems. 

\begin{figure}[h]
\begin{centering}
\subfloat[\label{fig:mbb_adap_result}]{\def\svgwidth{0.8\textwidth}     
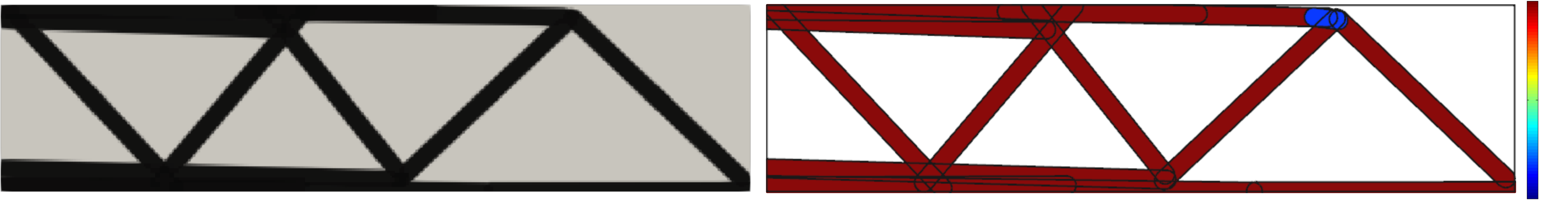

}
\par\end{centering}
\begin{centering}
\subfloat[\label{fig:mbb_full_result}]{\def\svgwidth{0.8\textwidth}     
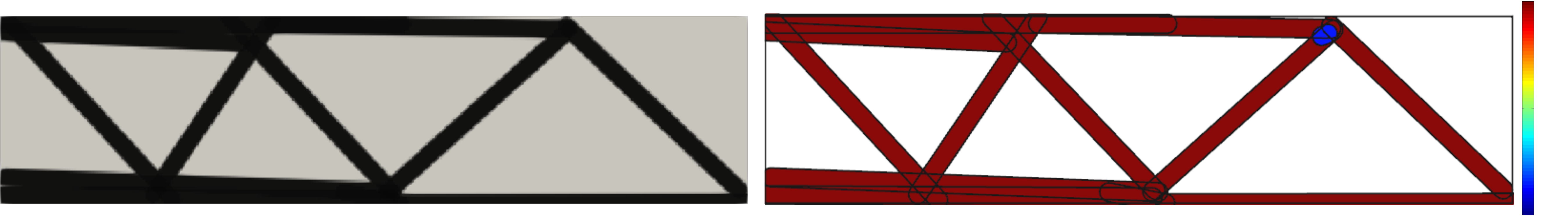

}
\par\end{centering}
\caption{(a) Optimal design obtained using (a) AMR; (b) full-resolution mesh.
Color denotes the penalized size variable $\alpha^{s}$.}
\end{figure}

Since the bars in the design are directly represented by the design
parameters, we can easily translate the optimal design into a CAD
model. The composite density field and its corresponding CAD model
for the optimal designs obtained using AMR and full-resolution meshes
are shown in Figs.\ \ref{fig:mbb_adap_result} and \ref{fig:mbb_full_result}
respectively. 

\begin{table}[h]
\begin{centering}
\begin{tabular}{lccc}
\hline 
Task & AMR (s) & Full resolution (s)  & AMR \% Improvement\tabularnewline
\hline 
Mesh refinement & 0.913 & 0 & -\tabularnewline
Geometry projection & 0.254 & 0.514 & 51\%\tabularnewline
Finite element assembly & 0.055 & 0.116 & 53\%\tabularnewline
Finite element linear solution & 0.114 & 0.784 & 85\%\tabularnewline
Responses calculation & 0.001 & 0.001 & 0\%\tabularnewline
Sensitivities calculation & 1.040 & 1.560 & 33\%\tabularnewline
\hline 
Total time & 2.337 & 2.975 & 21\%\tabularnewline
\hline 
\end{tabular}
\par\end{centering}
\caption{Time breakdown for the last iteration of the MBB problem.\label{tab:time}}
\end{table}

\subsection{Stress-constrained 2-dimensional L-bracket design }

\begin{figure}[h]
\center
\def\svgwidth{0.35\textwidth}     
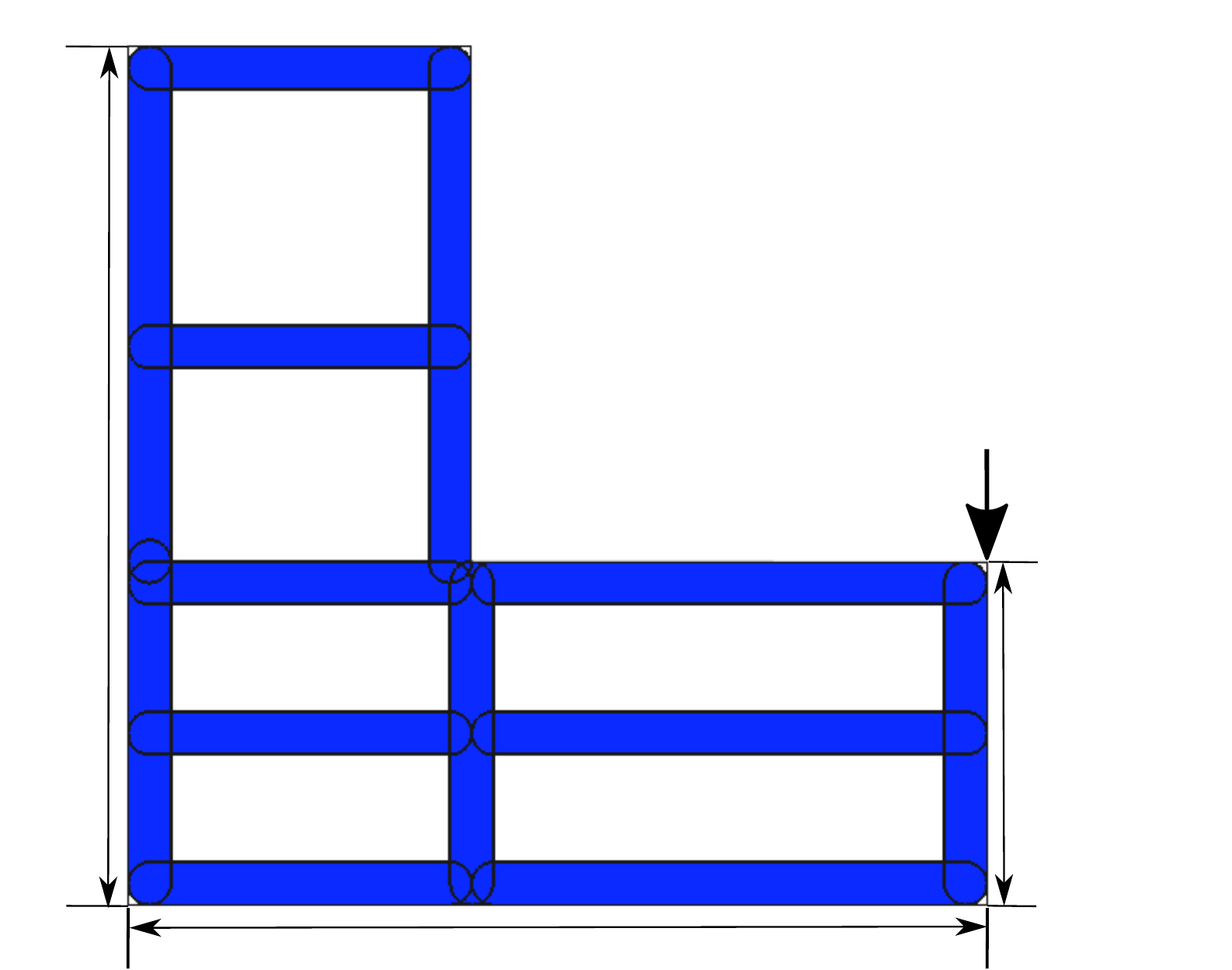\caption{Initial design, geometry, loads and boundary conditions for the 2-dimensional
L-bracket design problem. Color denotes the penalized size variable
$\alpha^{s}$. \label{fig:lshape-init}}

\end{figure}

In this example, we demonstrate the proposed AMR method in the context
of stress-based optimization. We consider the classical L-bracket
design problem that minimizes volume fraction with a constraint on
the highest von Mises stress in the structure. The optimization problem
is described in Eqs.\ \ref{eq:vol_obj} to \ref{eq:dv_bounds-2}.
The design envelope, loading, boundary conditions and initial design
are shown in Fig.\ \ref{fig:lshape-init}. The load $F$ has a magnitude
of 3 and is distributed on several elements along the vertical edge
to avoid the artificial stress concentration arising from a single
concentrated load. A maximum allowable stress constraint $\sigma^{\star}=2.4$
is imposed. We also impose a tight move limit $m=0.015$, to dampen
detrimental design changes arising from the nonlinearity of the problem
(cf., \cite{zhang2017stress}). Similarly to the previous example,
we perform the optimization on an adaptively refined mesh and a full-resolution
mesh for comparison. For AMR, the coarsest mesh level $L=0$ is a
uniform mesh with element size $h_{c}=2$. The target refinement level
$N_{rl}$ is set to 1 so that the finest mesh element has size of
$h=1$ after AMR. In order to have a consistent loading for both optimizations,
elements along the loading edge have an element size $h=1$ throughout
the optimization when using AMR. The full-resolution fine mesh has
a uniform mesh size $h=1$ with a total of 6400 elements. 

\begin{figure}[h]
\begin{centering}
\subfloat[\label{fig:2d-lshape-adap}]{\def\svgwidth{0.8\textwidth}     
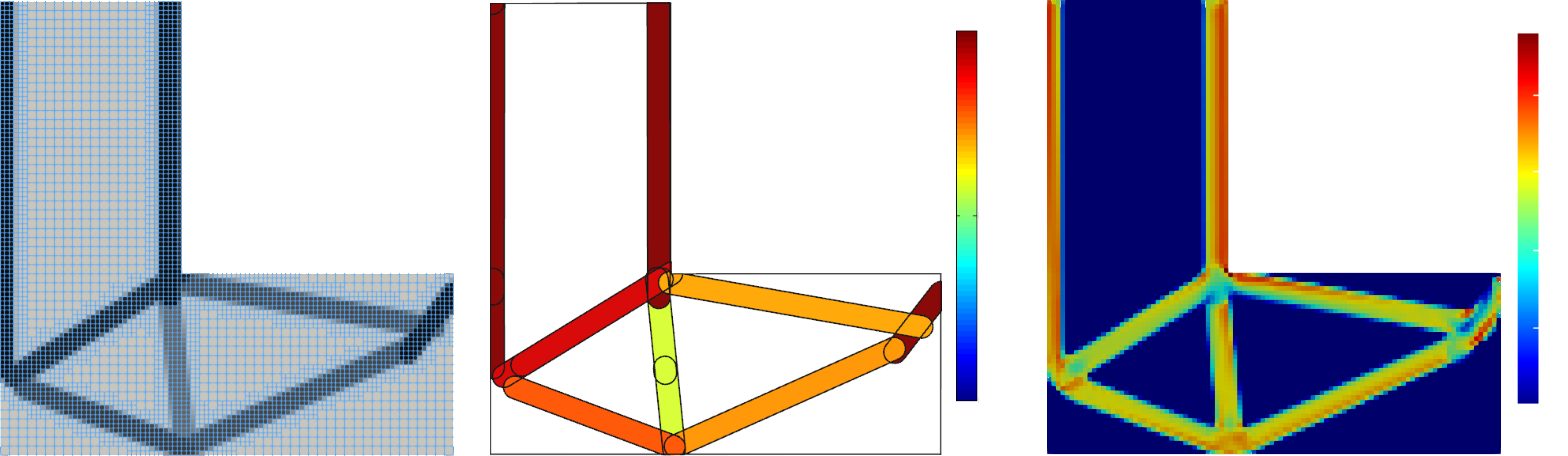

}
\par\end{centering}
\begin{centering}
\subfloat[\label{fig:2d-lshape-full}]{\def\svgwidth{0.8\textwidth}     
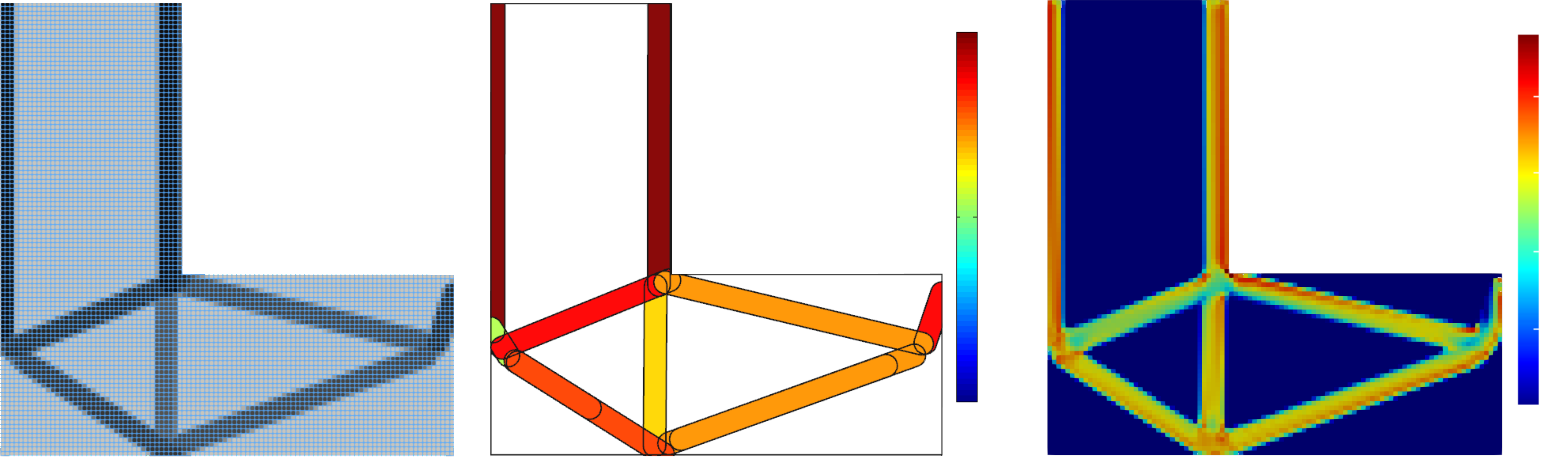

}
\par\end{centering}
\caption{(a) Optimal design and its relaxed stress distribution obtained using
(a) AMR, with $v_{f}=0.229$; and (b) a full-resolution mesh, with
$v_{f}=0.227$. \label{fig:2d-lshape-stress}}
\end{figure}

The optimization run with AMR takes 97 iterations to converge to a
design with $v_{f}=0.229$ as shown in Fig.\ \ref{fig:2d-lshape-adap}.
It takes 120 seconds to converge, and the final adaptively refined
mesh has 4006 elements, which roughly corresponds to two thirds the
number of elements of the full-resolution mesh. Fig.\ \ref{fig:2d-lshape-full}
shows for comparison the optimal design obtained with the full-resolution
mesh, which takes 117 iterations to converge to a design with $v_{f}=0.227$
in 159 seconds. We observe that both approaches obtain similar designs
with similar stress distributions. This example demonstrates that
the proposed AMR strategy can be applied to stress-based optimization.
We note that the recently published method in \cite{Troya2018} investigates
an AMR strategy to improve the accuracy of the stress computation
in density-based topology optimization. However, such a consideration
is out of the scope of this paper, and our goal is simply to reduce
the number of elements for large-scale problems.


\subsection{Stress-constrained 3-dimensional L-bracket design}

\begin{figure}[h]
\begin{centering}
\subfloat[\label{fig:3d-lshape-init}]{\def\svgwidth{0.28\textwidth}     
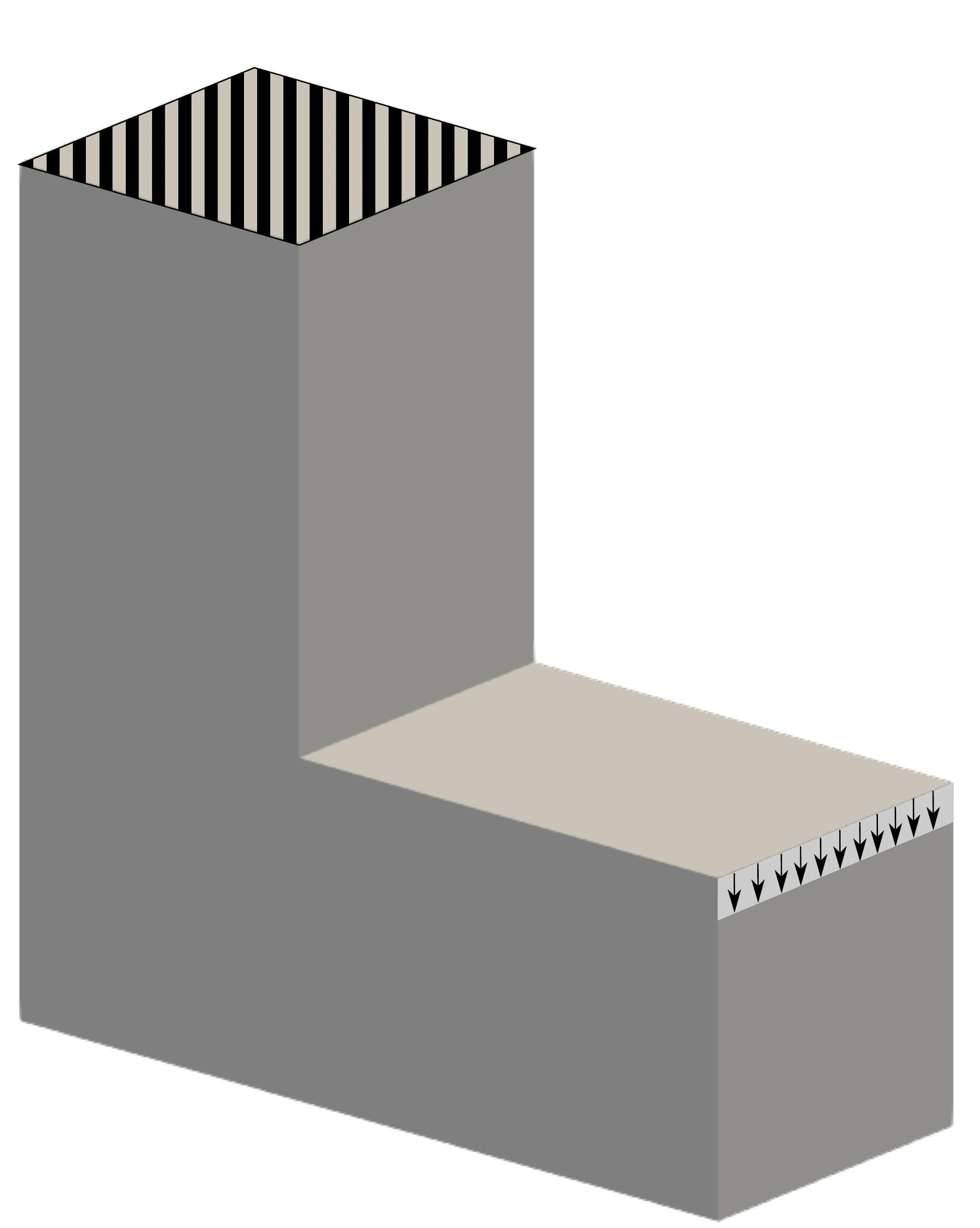

}\qquad{}\subfloat[]{\def\svgwidth{0.34\textwidth}     
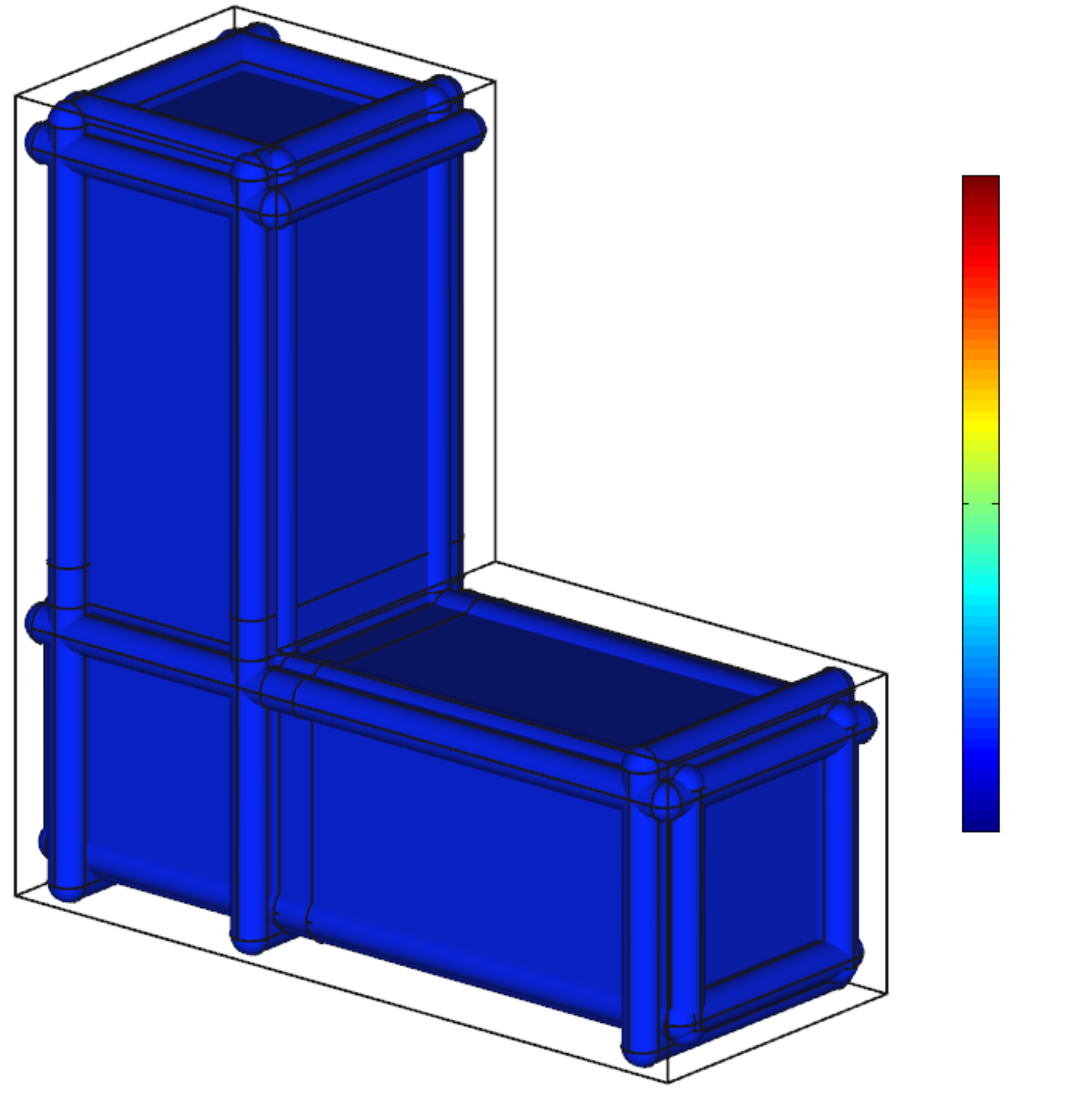

}\caption{3-dimensional L-bracket problem definition: (a) design space, loading
and boundary conditions; (b) initial design made of 16 plates. Color
scale indicates the penalized size variable $\alpha^{s}$.}
\par\end{centering}
\end{figure}

We now show another example of stress-based topology optimization
using AMR, this time on a 3-dimensional problem. We generate the design
envelope for this problem by extruding the 2-dimensional L-bracket
of the previous example by 40 units. The loading and boundary conditions
are shown in Fig.\ \ref{fig:3d-lshape-init}. A distributed surface
load is applied over the shaded area near the tip of the beam, as
depicted in Fig.\ \ref{fig:3d-lshape-init}, to avoid artificial
stress concentrations. The optimization consists of minimizing the
volume fraction $v_{f}$ with a maximum allowable stress constraint
$\sigma^{\star}=0.15$. The coarsest mesh level $L=0$ corresponds
to a uniform mesh with element size $h_{c}=4$, while the target refinement
level $N_{rl}$ is set to 2 so that the finest mesh element has size
of $h=1$ after AMR. As in the previous example, elements on the loading
surface are fixed and have an element size of $h=1$ throughout the
optimization when using AMR in order to have a consistent loading.
We employ the geometric multigrid preconditioner described in Sec.\ 
\ref{sec:Optimization-problems-and} for the solution of the analysis.
As in the previous examples, we perform two optimizations using AMR
and full-resolution meshes for comparison. The full-resolution mesh
is a uniform mesh with element size of $h=1$ and a total of 256,000
elements. For this example, we use four machines with 24 Xeon E5-2690
v3 2.60GHz cores each, and we employ one thread per core, for a total
of 96 threads. 

\begin{figure}[h]
\subfloat[]{\def\svgwidth{0.9\textwidth}     
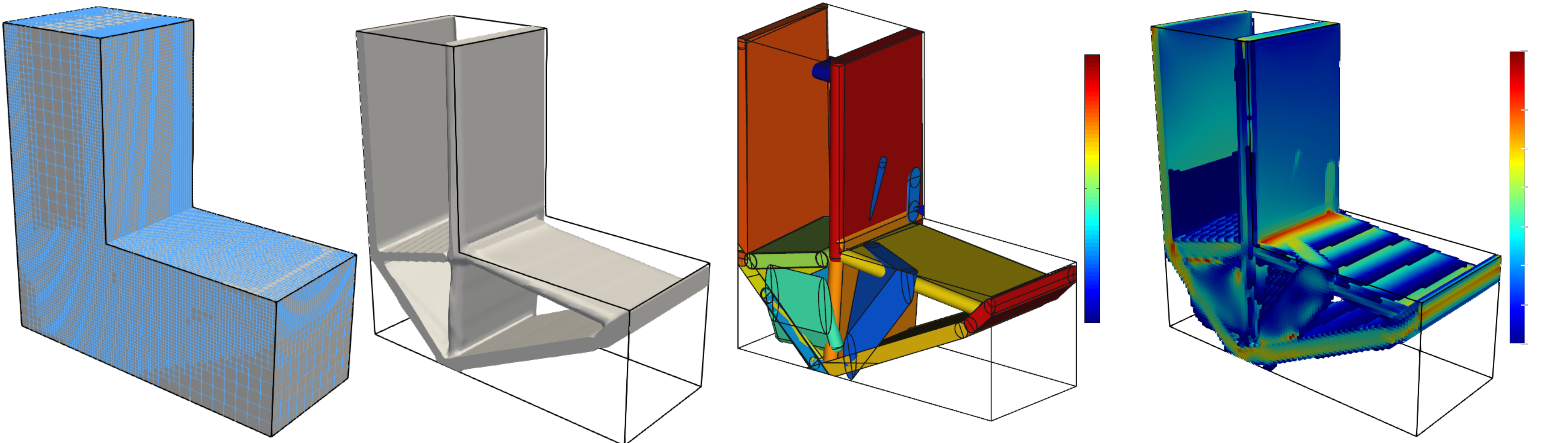

}

\subfloat[]{\def\svgwidth{0.9\textwidth}     
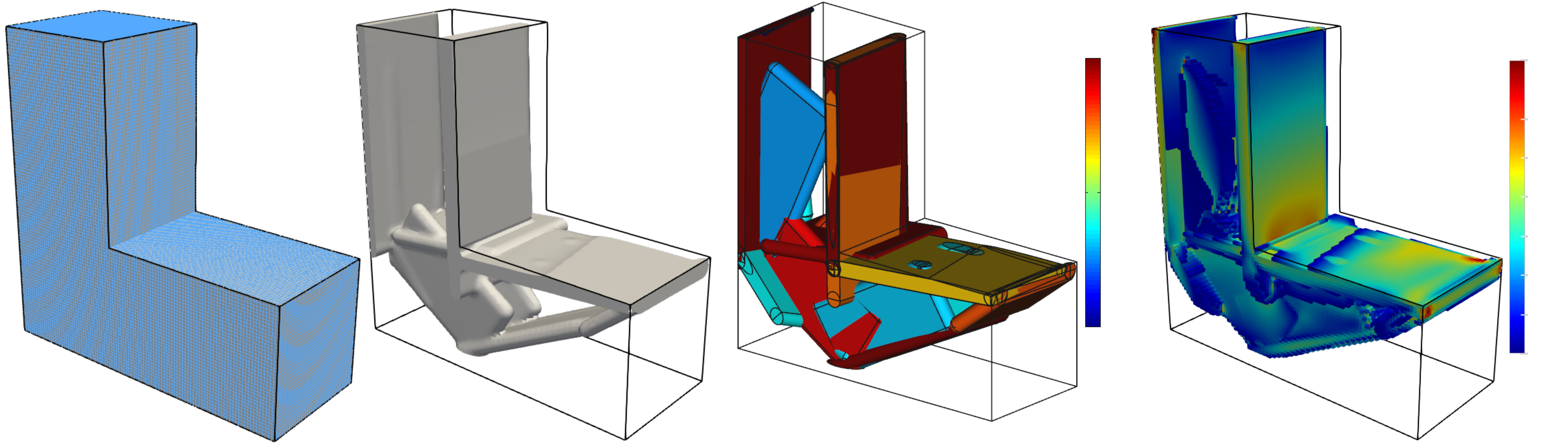

}

\caption{Results obtained for optimizations with (a) AMR, with $v_{f}=0.215$;
and (b) full-resolution mesh, with $v_{f}=0.212$. From left to right:
mesh for last iteration, composite density iso-surface of $\tilde{\rho}=0.4$,
CAD representation with color indicating penalized size variable $\alpha^{s}$,
and von Mises stress. \label{fig:3d-lshape-result}}

\end{figure}

Fig.\ \ref{fig:3d-lshape-result} shows the finite element mesh,
the composite density iso-surface of $\tilde{\rho}=0.4$, the CAD
model and the stress distribution of the optimal designs for the two
optimizations. The optimization with AMR takes 237 iterations to converge
to a design with $v_{f}=0.215$ in 18,000 seconds (5 hours). The adaptive
mesh for the optimal design using AMR has 157,321 elements. The optimization
with the full-resolution mesh converges to a design with $v_{f}=0.212$
in 213 iterations and 3,840 seconds (1.07 hours). It may seem surprising
at first that the optimization with the full-resolution mesh is much
faster than the one using AMR despite having far more elements. The
reason for this is that the local mesh refinement in AMR leads to
a less efficient GMG preconditioner than the one obtained for the
full-resolution mesh, which is obtained by globally refining the coarsest
mesh $N_{rl}$ times. This is evident from the comparison of the timing
information shown in Tab.\ \ref{tab:time-3d-lshape}. The AMR leads
to time savings for most of the tasks except the finite element linear
solution and the sensitivity analysis, which, when using iterative
solvers, requires another full-fledged finite element solution. Therefore,
we note that for 3-dimensional problems that have a modest size, like
the one presented here, the savings from using AMR cannot compensate
the loss in the effectiveness of the GMG preconditioner when compared
to the full-resolution uniform mesh. However, as the problem size
increases, the AMR with GMG will eventually outperform the full-resolution
mesh with GMG simply because of the substantial reduction in the number
of elements. 

\begin{table}[h]
\begin{centering}
\begin{tabular}{lccc}
\hline 
Task & AMR (s) & Full resolution (s)  & AMR \% Improvement\tabularnewline
\hline 
Mesh refinement & 1.78 & 0 & -\tabularnewline
Geometry projection & 0.238 & 0.259 & 8\%\tabularnewline
Finite element assembly & 0.275 & 0.398 & 31\%\tabularnewline
GMG preparation & 0.891 & 1.045 & 15\%\tabularnewline
Finite element linear solution & 55.6 & 10.2 & -445\%\tabularnewline
Responses calculation & 0.184 & 0.267 & 31\%\tabularnewline
Sensitivities calculation & 32.9 & 10.9 & -202\%\tabularnewline
\hline 
Total time & 91.868 & 23.069 & -298\%\tabularnewline
\hline 
\end{tabular}
\par\end{centering}
\caption{Time breakdown for the last iteration of the 3-dimensional L-bracket
problem.\label{tab:time-3d-lshape}}
\end{table}

\subsection{Compliance minimization of 3-dimensional cantilever beam }

\begin{figure}[h]
\begin{centering}
\subfloat[\label{fig:cantilever-initial}]{\def\svgwidth{0.45\textwidth}     
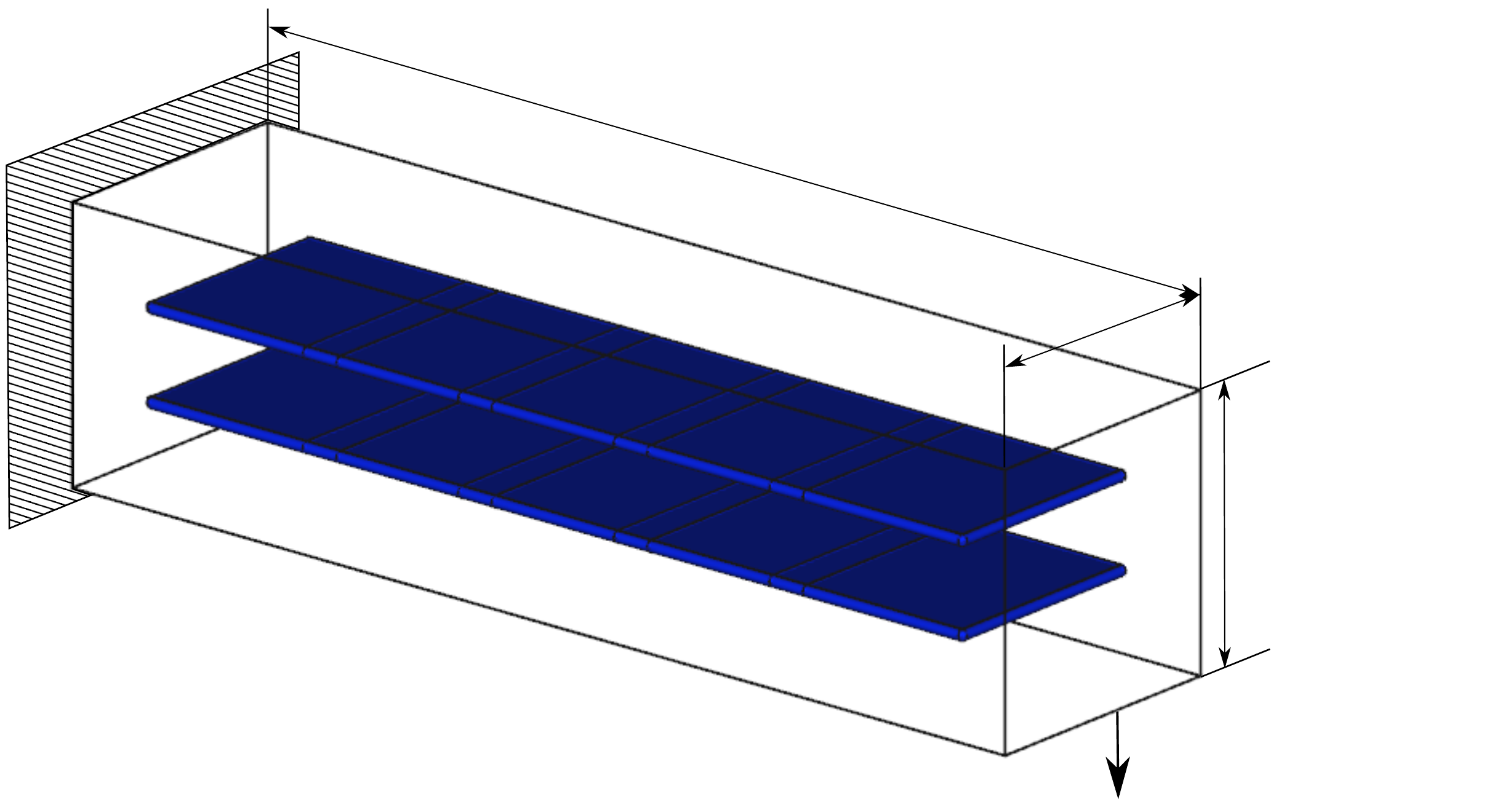

}\subfloat[\label{fig:cantilever-initial-mesh}]{\def\svgwidth{0.35\textwidth}     
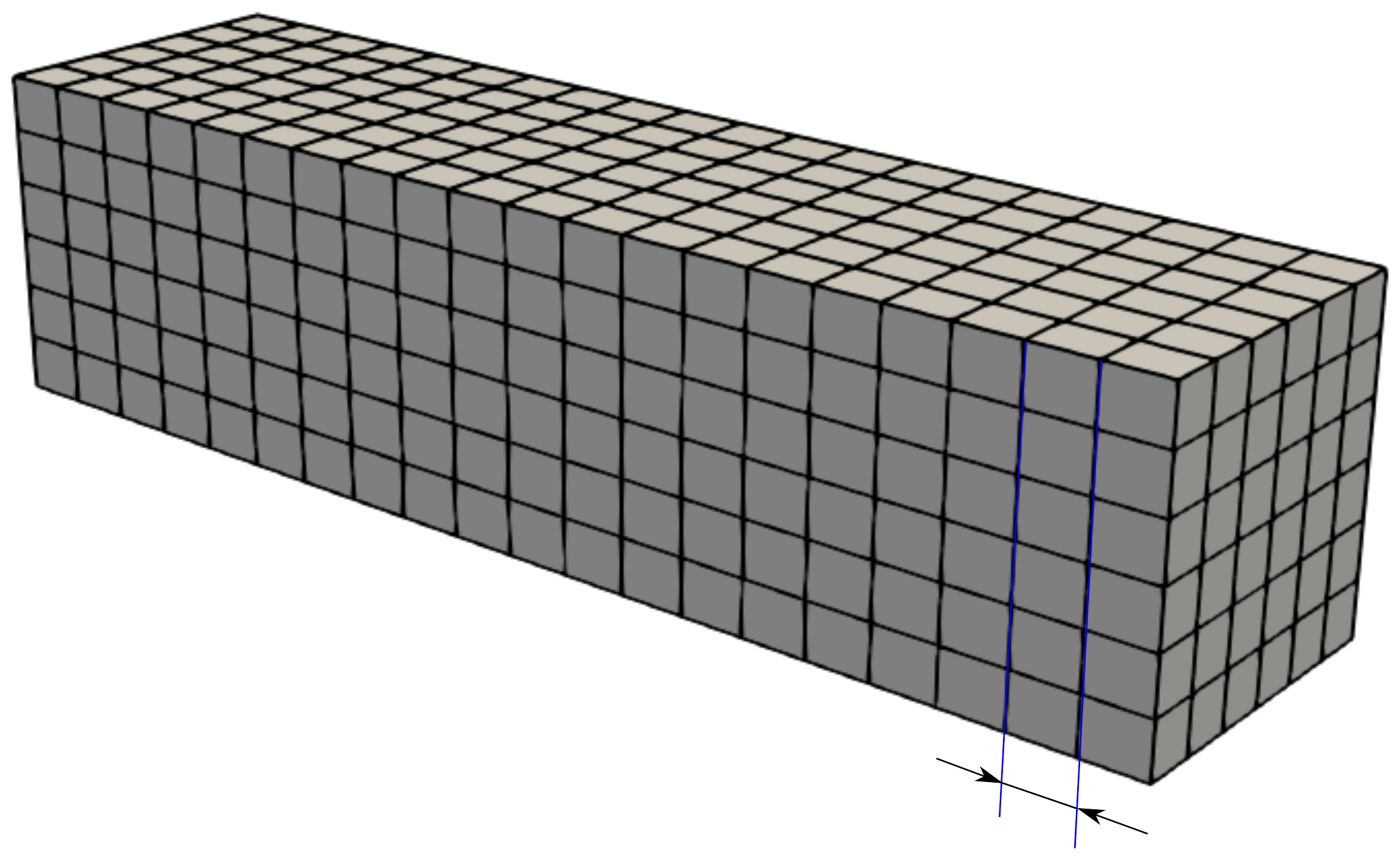

}
\par\end{centering}
\caption{(a) Initial design, geometry, loads and boundary conditions for the
3-dimensional cantilever beam design problem. Color denotes the penalized
size variable $\alpha^{s}$. (b) Initial coarse mesh.}

\end{figure}

In this last example, we consider the design of a 3-dimensional cantilever
beam to minimize its structural compliance subject to a relatively
low volume fraction constraint of $v^{*}=0.05$. The design envelope,
loading, boundary conditions and initial design are shown in Fig.\ \ref{fig:cantilever-initial}.
The initial design is made of eight plates. The applied load $F$
has a magnitude of 10. The MMA coefficient $M$ is set to 100 to better
satisfy the volume fraction constraint, and we use a move limit $m=0.05$.
We impose lower bounds $l_{\min}=2$ and $w_{\min}=2$ on the dimensions
of plates. The placement constraint introduced in \cite{zhang2016geometry}
is imposed for this example to ensure plates are entirely contained
within the design envelope and avoid impractical cuts across the plate
thickness (we refer the reader to the aforementioned publication for
details). All plates have a very small fixed thickness $t=0.2$ compared
to the dimensions of the design envelope. This high contrast between
the plates thickness and the dimensions of the design envelope is
common in the design of structures made of stock plates. For this
example, the coarsest mesh ($L=0$) has a uniform element size $h_{c}=1$
as shown in Fig.\ \ref{fig:cantilever-initial-mesh}. We target an
element size $h=0.03125$ that requires $N_{rl}=5$ levels of refinement.
If the mesh were globally refined five times to obtain the full-resolution
mesh as in previous examples, the refinement would result in 23,592,960
uniform mesh elements. However, by employing the proposed AMR strategy,
the number of elements reduces to 3,294,920 for the same problem (for
the initial design).

\begin{figure}[h]
\begin{centering}
\subfloat[]{\includegraphics[width=0.4\textwidth]{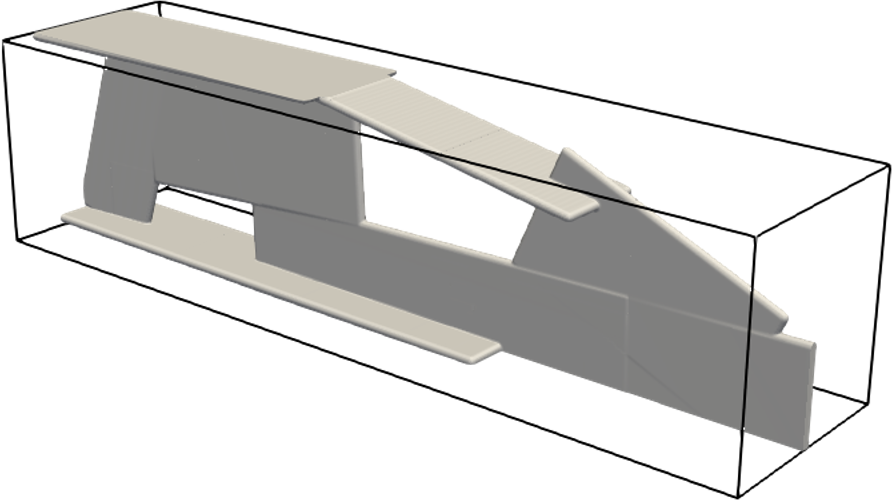}

}\subfloat[]{\def\svgwidth{0.45\textwidth}     
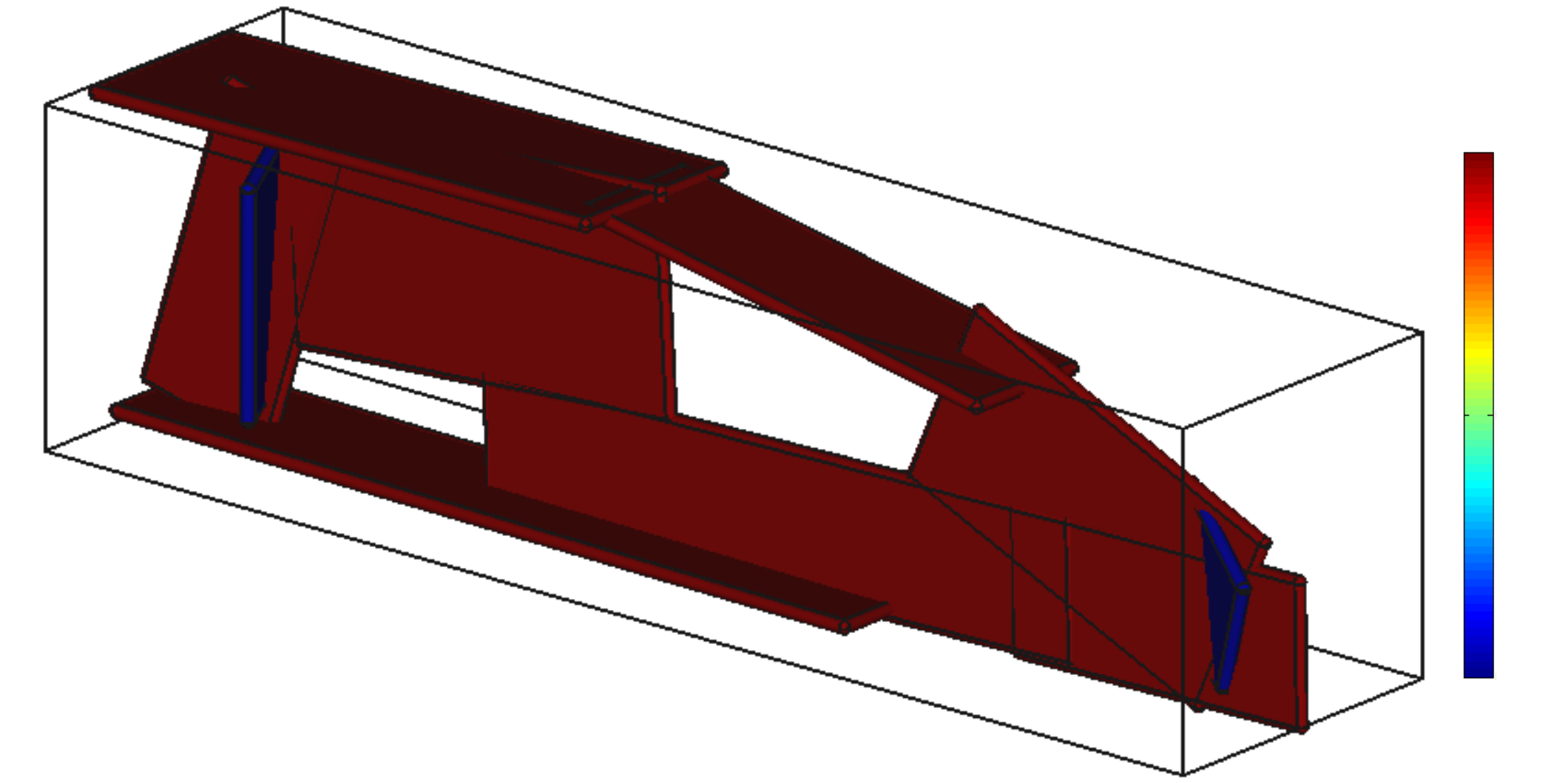

}
\par\end{centering}
\caption{Final design of 3-dimensional cantilever beam: (a) 0.4 iso-surface
of composite density, $C=2.71$; and (b) plate design, with colors
indicating the penalized size variable $\alpha^{s}$.\label{fig:can_final_design}}

\end{figure}

\begin{figure}[h]
\begin{centering}
\includegraphics[width=0.8\textwidth]{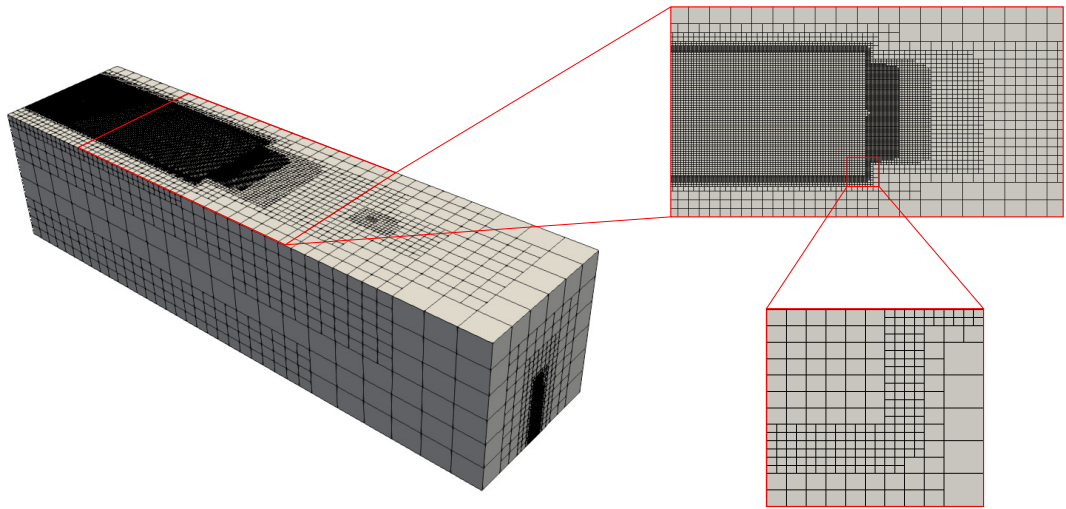}
\par\end{centering}
\caption{Adaptively refined mesh of the design in the last iteration of the
3-dimensional cantilever beam optimization using AMR.\label{fig:mesh-last-iter}}

\end{figure}

\begin{figure}[h]
\begin{centering}
\includegraphics[width=0.7\textwidth]{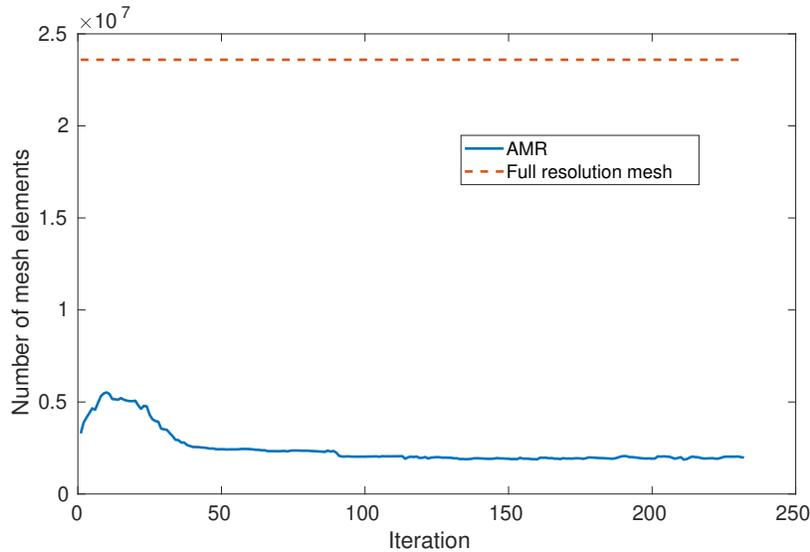}
\par\end{centering}
\caption{History of the number of elements using AMR for the 3-dimensional
cantilever design problem.\label{fig:nelem_history}}

\end{figure}

For this example, we use six machines with 24 Xeon E5-2690 v3 2.60GHz
cores each, and we employ one thread per core, for a total of 144
threads. We apply the AMR strategy to this example; Fig.\ \ref{fig:can_final_design}
shows the final design. The optimization takes 231 iteration in 665
minutes (11 h 5 min) to converge, hence an optimization iteration
spends 2.9 minutes on average. The adaptively refined mesh for the
last iteration is shown in Fig.\ \ref{fig:mesh-last-iter}. The total
number of elements for this last iteration is 1,976,155, which corresponds
to only 8.4\% of the full-resolution mesh. Fig.\ \ref{fig:nelem_history}
depicts the history of the number of mesh elements during the entire
optimization. For comparison, we also perform one design iteration
with the full-resolution mesh: we take the optimal design of Fig.\ \ref{fig:can_final_design}
obtained using AMR, refine all elements to the highest refinement
level to obtain the full-resolution mesh, and perform one design step.
If the optimization with the full-resolution mesh were to take the
same number of iterations as that with the AMR strategy, the total
optimization time would be approximately 42.5 hours (the actual time
will be different since the condition number of the system is different
for each design and therefore the number of iterations to convergence
of CG will vary). The timing information is summarized and compared
in Tab.\ \ref{tab:can-ex-time}. From this table, we can observe
that AMR leads to speedups in every task, resulting in an estimated
74\% time savings per optimization iteration. However, it is worth
noting that the savings obtained from AMR in the linear solution task
is not significant. The GMG preconditioner obtained from the full-resolution
mesh is so efficient that it can solve a 23 M elements mesh in approximately
seven minutes. However, the solution time with the AMR mesh is still
much shorter than with the full-resolution mesh. We note that in the
case of unstructured meshes, the efficiency of the GMG preconditioner
is limited, hence we expect a better speedup using AMR \cite{mavriplis1995multigrid}.

Unfortunately, there is no rule to determine a priori if performing
AMR is more efficient than using the full-resolution mesh, because
this depends on the problem; for instance, for a design region that
is not a cuboid, and/or for unstructured meshes as aforementioned,
computing the GMG preconditioner on the full-resolution mesh may become
more expensive, and so the AMR strategy may be advantageous for smaller
problems. The only certain way to determine the convenience of using
AMR is to perform a single finite element solution with both schemes
and compare the times. 

\begin{table}[h]
\begin{centering}
\begin{tabular}{lccc}
\hline 
Task & AMR (s) & Full resolution (s)  & AMR \% Improvement\tabularnewline
\hline 
Mesh refinement & 11 & 0 & -\tabularnewline
Geometry projection & 2.35 & 16.7 & 86\%\tabularnewline
Finite element assembly & 2.37 & 22.6 & 90\%\tabularnewline
GMG preparation & 6.16 & 69.22 & 91\%\tabularnewline
Finite element linear solution & 129 & 436 & 70\%\tabularnewline
Responses calculation & 0.003 & 0.019 & 84\%\tabularnewline
Sensitivities calculation & 4.58 & 51.7 & 91\%\tabularnewline
\hline 
Total time & 155.46 & 596.24 & 74\%\tabularnewline
\hline 
\end{tabular}
\par\end{centering}
\caption{Time breakdown for the last iteration of the 3-dimensional cantilever
beam problem.\label{tab:can-ex-time}}
\end{table}

\section{Conclusions}

The presented numerical experiments show that the proposed adaptive
mesh refinement strategy accommodates very well the topology optimization
with discrete geometric components using the geometry projection method.
The proposed strategy is independent of the shape of the geometric
components, and we applied it successfully to compliance minimization
and stress-constrained problems. Since the refinement indicator is
based on the geometry projection, our method can readily refine the
mesh at every optimization iteration, which leads to efficient designs.
As shown by our examples, it does not always make sense to use this
AMR strategy: for relatively small problems, or for large problems
where a highly structured mesh leads to a highly efficient GMG preconditioner,
the optimization with the full-resolution mesh may still outperform
the optimization with the adaptively refined mesh. However, for some
of the problems we are interested in, namely problems with very slender
members in relation to the dimensions of the design envelop and with
relatively low volume fractions, the required full-resolution meshes
are so large that the AMR strategy is bound to outperform the optimization
with the full-resolution mesh. In these cases, the AMR clearly achieves
our goal of substantially decreasing the computational burden. Additional
future work is required to incorporate in our refinement strategy
a refinement indicator based on solution error.

\section*{Acknowledgements}

Support from Caterpillar to conduct this work is gratefully acknowledged
by all authors. The last author would also like to acknowledge support
from the Office of Naval Research, USA, Grant N00014-17-1-2505.

\bibliography{AMR_arXiv}

\end{document}